%% file: Main.tex
\documentclass{CustomArticle}

\usepackage{algorithm}
\usepackage{algpseudocode}
\usepackage[style=trad-abbrv,citestyle=numeric-comp]{biblatex}
\usepackage{hyperref}
\usepackage[inline,shortlabels]{enumitem}
\usepackage{graphicx}
\usepackage{import}
\usepackage[capitalize,noabbrev,compress]{cleveref}

\usepackage{CommonMath}
\usepackage{ProjectMath}
\usepackage{TheoremParts}
\usepackage{LabelNamespaces}

\crefname{equation}{}{}

\newcounter{Halgorithmic}
\AtBeginEnvironment{algorithmic}{\stepcounter{Halgorithmic}}
\ExpandArgs{c}\newcommand{theHALG@line}{\arabic{Halgorithmic}.\arabic{ALG@line}}

\PaperGeneralInfo{
  title = {%
    Universal Gradient Methods\\%
    for Stochastic Convex Optimization%
  },
  date = {July 11, 2024}
}

\PaperAbstract{%
  We develop universal gradient methods for
  Stochastic Convex Optimization (SCO).
  Our algorithms automatically adapt not only to the oracle's noise
  but also to the Hölder smoothness of the objective function without
  a priori knowledge of the particular setting.
  The key ingredient is a novel strategy for adjusting step-size coefficients
  in the Stochastic Gradient Method (SGD).
  Unlike AdaGrad, which accumulates gradient norms, our
  Universal Gradient Method accumulates appropriate combinations of gradient-
  and iterate differences.
  The resulting algorithm has state-of-the-art worst-case convergence
  rate guarantees for the entire Hölder class including, in particular, both
  nonsmooth functions and those with Lipschitz continuous gradient.
  We also present the Universal Fast Gradient Method for SCO enjoying optimal
  efficiency estimates.
}

\AddPaperAuthor{
  name = {Anton Rodomanov},
  affiliation = {CISPA Helmholtz Center for Information Security},
  email = {anton.rodomanov@cispa.de}
}
\AddPaperAuthor{
  name = {Ali Kavis},
  affiliation = {%
    Institute for Foundations of Machine Learning (IFML),
    UT Austin%
  },
  email = {kavis@austin.utexas.edu}
}
\AddPaperAuthor{
  name = {Yongtao Wu},
  affiliation = {%
    Laboratory for Information and Inference Systems (LIONS),
    EPFL%
  },
  email = {yongtao.wu@epfl.ch}
}
\AddPaperAuthor{
  name = {Kimon Antonakopoulos},
  affiliation = {%
    Laboratory for Information and Inference Systems (LIONS),
    EPFL%
  },
  email = {kimon.antonakopoulos@epfl.ch}
}
\AddPaperAuthor{
  name = {Volkan Cevher},
  affiliation = {%
    Laboratory for Information and Inference Systems (LIONS),
    EPFL%
  },
  email = {volkan.cevher@epfl.ch}
}

\PaperKeywords{%
  convex optimization,
  stochastic optimization,
  gradient methods,
  universal methods,
  adaptive algorithms,
  efficiency estimates,
  Hölder smoothness
}

\begin{document}
  \PrintTitleAndAbstract

  \subimport{Introduction/}{Main}
  \subimport{Preliminaries/}{Main}
  \subimport{GradientMethod/}{Main}
  \subimport{StochasticGradientMethod/}{Main}
  \subimport{./}{StochasticFastGradientMethod}
  \subimport{Experiments/}{Main}
  \subimport{./}{Conclusion}

  \printbibliography

  \appendix
  \subimport{Appendix/}{Main}
\end{document}

%% file: Introduction/Main.tex
\section{Introduction}

\subimport{./}{Motivation}
\subimport{./}{Contributions}
\subimport{./}{RelatedWork}

%% file: Introduction/Motivation.tex
\paragraph{Motivation.}

The complexity of modern machine learning problems makes it difficult to
estimate their mathematical properties, let alone characterize them accurately.
The problems thus demand sophisticated solutions which are robust to possible
variations in the parameters.
One therefore needs algorithms that can work simultaneously under multiple
assumptions while implicitly adapting to the parameters of the problem.
The sheer scale of modern problems also raises efficiency
concerns, which paves the way for the stochastic methods leveraging
randomized computations.

In this paper, we study the convex optimization problem
\begin{equation}
  \label{eq:MainProblem-Informal}
  \min_{x \in \EffectiveDomain \psi}
  \bigl[ F(x) \DefinedEqual f(x) + \psi(x) \bigr],
\end{equation}
where $f$ is the main (difficult) part of the problem,
and $\psi$ is a \emph{simple} convex function (e.g., indicator of a set).
Furthermore, we assume that $f$ can be queried only via an unbiased
\emph{stochastic gradient oracle} with (unknown) variance~$\sigma^2$.

Optimization algorithms are typically designed for
a particular problem class and tailored to its properties.
Two standard classes are \emph{nonsmooth} ($f$ is Lipschitz continuous)
and \emph{smooth} ($f$ has Lipschitz gradient).
It is common for the problem class to dictate the selection of algorithm's
parameters to ensure the optimal convergence.

However, in practice, every specific problem typically belongs to multiple
problem classes at the same time, and it is usually very difficult
(if not impossible) to say in advance which particular class better fits our
problem.
To address this issue, we need \emph{universal methods} that can automatically
adjust to the ``correct'' problem class when applied to a concrete problem
instance given to them.

The important example of such algorithms is given by
Universal Gradient Methods (UGMs) of~\cite{nesterov2015universal}.
These methods are capable of solving the more general class of
\emph{Hölder-smooth} problems:
\[
  \Norm{\Gradient f(x) - \Gradient f(y)} \leq L_{\nu} \Norm{x - y}^{\nu},
  \quad
  \forall x, y \in \EffectiveDomain \psi,
\]
which continuously connects nonsmooth problems ($\nu = 0$)
with the smooth ones ($\nu = 1$).
To achieve universality, UGMs use a special line-search procedure
which automatically selects an appropriate step size for any possible Hölder
exponent and the corresponding Hölder constant, without knowing these
parameters.
As a result, the methods automatically adjust to the best possible problem
class.
However, UGMs require exact computations of gradients ($\sigma = 0$).

The extension of UGMs to stochastic optimization has been a challenging
\emph{open problem}.
The desired algorithms should automatically adjust not only to the Hölder
smoothness of the objective function, but also to the oracle's noise.

In this paper, we address this open problem and provide a solution to it.
We design line-search-free variants of UGMs which automatically adapt to:
\begin{enumerate*}[(\itshape i)]
  \item Hölder exponent~$\nu$,
  \item Hölder constant~$L_{\nu}$,
  \item variance of the stochastic oracle~$\sigma$,
\end{enumerate*}
without having the prior knowledge of neither the problem class nor the nature
of the gradient information.

%% file: Introduction/Contributions.tex
\paragraph{Contributions.}

We develop new Universal Gradient Methods (UGMs) for
problem~\eqref{eq:MainProblem-Informal}, which are robust to the stochastic
noise in gradient computations.
To achieve that, we assume the knowledge of a certain upper bound~$D$ on the
diameter of the feasible set~$\EffectiveDomain \psi$
(or, somewhat equivalently, the distance from the initial point to the solution),
which is a common assumption in a variety of other adaptive algorithms for
Stochastic Convex Optimization (SCO).

Our main contributions can be summarized as follows:
\begin{enumerate}
  \item
    We first rethink (in \cref{sec:GradientMethod}) the theoretical analysis
    of the line-search-based UGM for deterministic optimization, and identify
    a simple mechanism to remove the line search from this algorithm
    while retaining the same worst-case efficiency estimates.
    The key element is a \emph{novel strategy for adjusting step-size coefficients}
    based on the idea of balancing the two error terms appearing in the
    convergence analysis.
  \item
    We then show (in \cref{sec:GradientMethodForStochasticOptimization}) that
    our techniques can easily be extended to stochastic optimization problems.
    The only essential change that we need to make is to replace the
    Bregman distance for the objective function,
    appearing in the formula for the step-size, with the stochastic version
    of the symmetrized Bregman distance involving gradient- and iterate
    differences.
    The resulting Universal Stochastic Gradient Method requires at most
    $
      \BigO\bigl(
        \inf_{\nu \in \ClosedClosedInterval{0}{1}}
        [\frac{L_{\nu}}{\epsilon}]^{2 / (1 + \nu)} D^2
        +
        \frac{\sigma^2 D^2}{\epsilon^2}
      \bigr)
    $
    stochastic oracle calls to reach $\epsilon$-accuracy in terms of the
    expected function residual
    (\cref{th:ConvergenceRateOfStochasticGradientMethod}).

  \item
    Finally, we present
    (in \cref{sec:FastGradientMethodForStochasticOptimization})
    the Universal Stochastic Fast Gradient Method enjoying
    the worst-case optimal efficiency of
    $
      \BigO\bigl(
        \inf_{\nu \in \ClosedClosedInterval{0}{1}}
        [\frac{L_{\nu} D^{1 + \nu}}{\epsilon}]^{2 / (1 + 3 \nu)}
        +
        \frac{\sigma^2 D^2}{\epsilon^2}
      \bigr)
    $
    oracle calls (\cref{th:ConvergenceRateOfStochasticFastGradientMethod}).
\end{enumerate}

Note that all our methods are agnostic to the smoothness exponent~$\nu$,
smoothness constant~$L_{\nu}$ and the noise level~$\sigma$.
To our knowledge, this is the first work proposing algorithms with such
characteristics.

%% file: Introduction/RelatedWork.tex
\paragraph{Related work.}

Pioneered by the AdaGrad algorithm~\textcite{duchi2011adaptive,
mcmahan2010adaptive}, \emph{adaptive methods} have been at the forefront of
training machine learning models.
AdaGrad accumulates the sequence of observed gradient norms to construct a
decreasing step size.
This construction enables data-adaptive regret bounds and has many useful
properties.
Following the success of the AdaGrad, several methods have been
proposed~\parencite{%
  kingma2014adam,%
  tieleman2012lecture,%
  rakhlin2013optimization,%
  reddi2018convergence%
}.
\textcite{levy2018online} proposed the first accelerated algorithm with
data-adaptive step-size without the knowledge of Lipschitz constant
and the variance bound.
They prove convergence results for nonsmooth and smooth objectives in the
presence of stochastic noise.
These results are further refined and extended by~\textcite{%
  kavis2019unixgrad,%
  joulani2020simpler,%
  ene2021adaptive%
}.
Despite the significant interest in the adaptation to smoothness and noise,
existing methods are not known to handle Hölder-smooth objectives.

Another popular type of adaptive methods is known as \emph{parameter-free}.
This direction is very interesting but somewhat orthogonal to ours.
Parameter-free algorithms have been studied for over a decade in online
learning~\parencite{%
  mcmahan2012noregret,%
  orabona2014simultaneous,%
  cutkosky2017online,%
  cutkosky2018blackbox,%
  jacobsen2023unconstrained,%
  mhammedi2020lipschitz%
}.
They are usually endowed with appropriate mechanisms to achieve efficiency
bounds that are almost insensitive (typically, with logarithmic dependency)
to the error of estimating certain problem parameters, such as the diameter of
the feasible set~\parencite{%
  carmon2022making,%
  ivgi2023dog,%
  defazio2023learning,%
  khaled2023dowg,%
  mishchenko2023prodigy%
}.
However, these methods typically consider the extreme cases of the Hölder class.

Within the context of online learning, there exists an independent notion
of universality such that the algorithms adapt unknowingly to the degrees and
types of convexity.
The goal is designing algorithms that achieve, up to logarithmic factors,
optimal regret bounds simultaneously for convex, strongly convex and
exponentially concave functions
\cite{erven2016metagrad,wang2020adaptivity,zhang2022simple,yan2023universal}.
The associated design and proof techniques are not transferable to our setup as
we focus on the degree of smoothness while the aforementioned works study
degrees of convexity.

The first UGM for deterministic optimization, including the Fast UGM with
optimal worst-case oracle complexity, was proposed
in~\cite{nesterov2015universal}.
The corresponding methods achieve the adaptation to Hölder smoothness by
the means of line search but must set the target accuracy a priori.
A possible extension of these algorithms to stochastic optimization was
considered in~\cite{gasnikov2018universal}.
They proposed an accelerated
gradient method for stochastic optimization problems, which adapts to the Hölder
characteristics of the objective using line search combined with mini-batching.
However, this method additionally relies on the knowledge of the oracle's
variance to correctly set up the size of the mini-batch at each iteration,
and therefore cannot be considered adaptable to the noise level.

More recently, \textcite{li2023simple} studied the same problem but with the
deterministic oracle, and designed a line-search-free universal method that
estimates local smoothness in the sense
of~\textcite{malitsky2020adaptive,malitsky2023adaptive}.
Their step-size formula shares some similarities to ours, and does not
require any (artificial) bounds on the diameter of the feasible set.
However, they only consider exact gradient computations, and it is unknown
whether their construction can be extended to stochastic problems.
On a related note, \textcite{orabona2023normalized} showed that, in the
deterministic case, both AdaGrad and the normalized gradient
method~\parencite[Section 3.2.3]{nesterov2018lectures} automatically adapt to
Hölder smoothness.
Although the corresponding proof for AdaGrad could be extended to the stochastic
setting at the expense of additional assumptions, the same type of argument
cannot be trivially applied to the accelerated method.

%% file: Preliminaries/Main.tex
\section{Preliminaries}

\subimport{./}{Notation}
\subimport{./}{ProblemSetting}

%% file: Preliminaries/Notation.tex
\subsection{Notation}

In this text, we work in the space $\RealField^n$ equipped with the standard
inner product $\DualPairing{\cdot}{\cdot}$ and the certain Euclidean norm:
\begin{equation}
  \label{eq:EuclideanNorm}
  \Norm{x} \DefinedEqual \DualPairing{B x}{x}^{1 / 2},
  \qquad
  x \in \RealField^n,
\end{equation}
where $B \in \SymmetricPdMatrices{n}$ is a sufficiently simple symmetric
positive definite matrix (e.g., the identity or a diagonal one).
The corresponding dual norm is defined in the standard way:
\begin{equation}
  \label{eq:DualNorm}
  \DualNorm{s}
  \DefinedEqual
  \max_{\Norm{x} = 1} \DualPairing{s}{x}
  =
  \DualPairing{s}{B^{-1} s}^{1 / 2},
  \qquad
  s \in \RealField^n.
\end{equation}
Thus, for any $s, x \in \RealField^n$, we have the Cauchy--Schwarz
inequality $\Abs{\DualPairing{s}{x}} \leq \DualNorm{s} \Norm{x}$.

For a convex function $\Map{f}{\RealField^n}{\RealFieldPlusInfty}$, by
$\EffectiveDomain f \DefinedEqual \SetBuilder{x \in \RealField^n}{f(x) <
+\infty}$, we denote its \emph{effective domain}.
The subdifferential of $f$ at a point $x \in \EffectiveDomain f$ is denoted
by~$\Subdifferential f(x)$.
For any two points $x, y \in \EffectiveDomain f$ and any $g \in \Subdifferential
f(x)$, we define the \emph{Bregman distance} generated by~$f$ as
\begin{equation}
  \label{eq:BregmanDistance}
  \BregmanDistanceWithSubgradient{f}{g}(x, y)
  \DefinedEqual
  f(y) - f(x) - \DualPairing{g}{y - x}
  \quad (\geq 0),
\end{equation}
In the case when there is no ambiguity with the subgradient~$g$,
we use a simpler notation $\BregmanDistance{f}(x, y)$.

For any $t \in \RealField$, by $\PositivePart{t} \DefinedEqual \max \Set{t, 0}$,
we denote its positive part.
For random variables $X$ and $\xi$, by $\Expectation_{\xi}[X]$ and
$\Expectation[X]$, we denote the expectation of $X$ w.r.t.\ to $\xi$, and the
full expectation of $X$, respectively.

%% file: Preliminaries/ProblemSetting.tex
\subsection{Problem Setting}

In this paper, we study the following optimization problem:
\begin{equation}
    \label{eq:MainProblem}
    F^*
    \DefinedEqual
    \min_{x \in \EffectiveDomain \psi}
    \bigl[ F(x) \DefinedEqual f(x) + \psi(x) \bigr],
\end{equation}
where $\Map{\psi}{\RealField^n}{\RealFieldPlusInfty}$ is a sufficiently
\emph{simple} proper closed convex function, and
$\Map{f}{\RealField^n}{\RealFieldPlusInfty}$ is a closed convex function
which is finite and subdifferentiable over an open set containing
$\EffectiveDomain \psi$.

Our main assumption on problem~\eqref{eq:MainProblem} is the
\emph{boundedness of the feasible set}~$\EffectiveDomain \psi$.

\begin{assumption}
  \label{as:BoundedFeasibleSet}
  For problem~\eqref{eq:MainProblem}, there is $D > 0$ such that
  $\Norm{x - y} \leq D$ for all $x, y \in \EffectiveDomain \psi$.
\end{assumption}

In what follows, we assume that the diameter~$D$ is known.
This will be the only parameter in our methods.
Note that \cref{as:BoundedFeasibleSet} guarantees that the
problem~\eqref{eq:MainProblem} has a solution (since the objective function~$F$
is proper and closed).
An important example that satisfies this assumption is when $\psi$ is the
indicator function of a certain compact set $Q \subseteq \RealField^n$:
$\psi(x) \DefinedEqual 0$ if $x \in Q$, and $\psi(x) \DefinedEqual +\infty$ if $x \notin Q$.

By calling $\psi$ simple, we mean the following standard assumption:
for any $c \in \RealField^n$, $\bar{x} \in \EffectiveDomain \psi$, and
$H \geq 0$, we can efficiently compute a solution to the following subproblem:
$
    \min_{x \in \EffectiveDomain \psi} \{
      \DualPairing{c}{x}
      +
      \frac{H}{2} \Norm{x - \bar{x}}^2
      +
      \psi(x)
    \}
$.
For instance, when $\psi$ is the indicator function of a compact set~$Q$,
it corresponds to finding a Euclidean projection onto~$Q$ (or minimizing a
linear function over~$Q$ if $H = 0$; in this case, we allow for an arbitrary
solution of the subproblem).

To characterize the smoothness of~$f$ in
problem~\eqref{eq:MainProblem}, let us introduce the \emph{Hölder constant}
for each $\nu \in \ClosedClosedInterval{0}{1}$:
\begin{equation}
  \label{eq:HolderConstant}
  L_{\nu}
  \DefinedEqual
  \sup_{
    \substack{
      x, y \in \EffectiveDomain \psi, \
      x \neq y,
      \\
      g(x) \in \Subdifferential f(x), \
      g(y) \in \Subdifferential f(y)
    }
  }
  \frac{\DualNorm{g(x) - g(y)}}{\Norm{x - y}^{\nu}}.
\end{equation}
Of course, for certain values of the
exponent~$\nu \in \ClosedClosedInterval{0}{1}$, it may happen that $L_{\nu} = +\infty$.
However, we assume that there exists (at least one) exponent for which the
corresponding Hölder constant is finite.

\begin{assumption}
  \label{as:AtLeastOneHolderConstantIsFinite}
  For problem~\eqref{eq:MainProblem} and $L_{\nu}$ given by
  \cref{eq:HolderConstant}, there exists $\nu \in \ClosedClosedInterval{0}{1}$
  such that $L_{\nu} < +\infty$.
\end{assumption}

The case $L_0 < +\infty$ corresponds to the situation when~$f$ has
\emph{bounded variation of subgradients}: for all $x, y \in \EffectiveDomain \psi$,
and all $g(x) \in \Subdifferential f(x)$, $g(y) \in \Subdifferential f(y)$,
it holds that
$
  \DualNorm{g(x) - g(y)} \leq L_0.
$
If $f$ has bounded subgradients over $\EffectiveDomain \psi$,
i.e., there exists $L_0' \geq 0$ such that $\DualNorm{g(x)} \leq L_0'$ for all $x
\in \EffectiveDomain \psi$ and all $g(x) \in \Subdifferential f(x)$, then $L_0
\leq 2 L_0'$.
On the other hand, if $L_{\nu} < +\infty$ for some $\nu \in
\OpenClosedInterval{0}{1}$, then $f$ is actually differentiable over
$\EffectiveDomain \psi$,  and, for all $x, y \in
\EffectiveDomain \psi$, it holds
$
  \DualNorm{\Gradient f(x) - \Gradient f(y)}
  \leq
  L_{\nu} \Norm{x - y}^{\nu}.
$
The case $L_1 < +\infty$ corresponds to the \emph{Lipschitz gradient}.

One simple example of the convex function with Hölder (sub)gradients is
the $p$-th power of the $\ell_p$-norm of the residual for the system of
linear equations,
$
 f(x)
 =
 \frac{1}{m} \sum_{i = 1}^m \Abs{\langle a_i, x \rangle - b_i}^p
$
with $a_i \in \RealField^n$, $b_i \in \RealField$,
$p \in \ClosedClosedInterval{1}{2}$,
which generalizes the classical least-squares loss (corresponding to $p = 2$);
this function is Hölder smooth with $\nu = p - 1$ but not Lipschitz smooth
(unless $p = 2$).
Another simple example is the similar residual but for linear inequalities,
$
  f(x)
  =
  \frac{1}{m} \sum_{i = 1}^m \PositivePart{\langle a_i, x \rangle - b_i}^p
$,
which is the smooth counterpart of the classical loss function used by
the Support Vector Machines (SVMs).
More generally, there is a duality relationship between Hölder smoothness
and uniform convexity:
if $f_*$ is a uniformly convex function\footnote{%
  This means that
  $\InnerProduct{f_*'(x) - f_*'(y)}{x - y} \geq \sigma_q \Norm{x - y}^q$
  for all $x, y$ and all $f_*'(x) \in \Subdifferential f_*(x)$,
  $f_*'(y) \in \Subdifferential f_*(y)$.
}
of degree~$q \geq 2$ with parameter~$\sigma_q > 0$, then its Fenchel
dual $f$ is Hölder smooth with $\nu = \frac{1}{q - 1}$ and
$L_{\nu} \leq \bigl( \frac{1}{\sigma_q} \bigr)^{\frac{1}{q - 1}}$
(see, e.g., Lemma~1 in~\cite{nesterov2015universal}), and vice versa;
in particular\footnote{%
  Here we use the standard fact that $\frac{1}{q} \DualNorm{\cdot}^q$
  is a uniformly convex function of degree~$q$ with
  parameter~$\frac{1}{2^{q - 2}}$,
  (see, e.g., Lemma~4 in~\cite{nesterov2008accelerating}).
},
for any convex function $f_*$, the function
$
 f(x)
 =
 \max_{s \in \EffectiveDomain f_*} [
  \langle s, x \rangle - f_*(s) - \frac{\sigma_q}{q} \DualNorm{s}^q
 ]
$
is Hölder smooth with $\nu = \frac{1}{q - 1}$ and
$L_{\nu} \leq \bigl( \frac{2^{q - 2}}{\sigma_q} \bigr)^{\frac{1}{q - 1}}$.

It is not difficult to see from \cref{eq:HolderConstant} that, under
\cref{as:BoundedFeasibleSet}, for any $0 \leq \nu_1 \leq \nu_2 \leq 1$, we have
the following monotonicity relation:
$
  L_{\nu_1} D^{\nu_1}
  \leq
  L_{\nu_2} D^{\nu_2}.
$
(This is the consequence of the fact that $\tau^p$ is increasing in $p > 0$ for
any fixed $\tau = \frac{D}{\Norm{x - y}} \geq 1$.)
In particular, if $L_{\nu'} < +\infty$ for some $\nu' \in
\ClosedClosedInterval{0}{1}$, then $L_{\nu} < +\infty$ for all $\nu \in
\ClosedClosedInterval{0}{\nu'}$.

One standard and important consequence of \cref{eq:HolderConstant} is that, for
any $\nu \in \ClosedClosedInterval{0}{1}$ (such that $L_{\nu} < +\infty$), and
all $x, y \in \EffectiveDomain \psi$ and all $g \in \Subdifferential f(x)$, we
have the following upper bound on the Bregman distance of the function~$f$:
\begin{equation}
    \label{eq:UpperBoundOnBregmanDistanceViaHolderConstant}
    \BregmanDistanceWithSubgradient{f}{g}(x, y)
    \leq
    \frac{L_{\nu}}{1 + \nu} \Norm{x - y}^{1 + \nu}.
\end{equation}

Our goal in this paper is to present numerical methods for
solving~\eqref{eq:MainProblem} that are \emph{universal}: they can
automatically adapt to the actual level of smoothness of the function~$f$
without knowing neither the Hölder exponent~$\nu$, nor the corresponding Hölder
constant~$L_{\nu}$.

%% file: GradientMethod/Main.tex
\section{Universal Line-Search-Free Gradient Method}
\label{sec:GradientMethod}

\subimport{MainIdea/}{Main}
\subimport{./}{TheMethod}

%% file: GradientMethod/MainIdea/Main.tex
\subsection{Main Idea}
\label{sec:MainIdeaForGradientMethod}

\UsingNamespace{MainIdeaForGradientMethod}

To explain the main idea behind our construction of adaptive step-size
coefficients, let us consider the usual (Composite) Gradient Method for solving
problem~\eqref{eq:MainProblem}:
\begin{equation}
  \LocalLabel{eq:GradientMethod}
  x_{k + 1}
  =
  \argmin_{x \in \EffectiveDomain \psi} \Bigl\{
    \InnerProduct{f'(x_k)}{x}
    +
    \psi(x)
    +
    \frac{H_k}{2} \Norm{x - x_k}^2
  \Bigr\},
\end{equation}
assuming we can compute the exact (sub)gradient~$f'(x_k) \in \Subdifferential
f(x_k)$ at each iteration~$k \geq 0$ (i.e., the oracle is deterministic).
The question is how to choose the step-size coefficients~$H_k$ at each iteration
to ensure that the algorithm properly works for any possible Hölder
exponent~$\nu$ and the corresponding coefficient~$L_{\nu}$ without explicitly
using these constants in the method.

The standard convergence analysis of
method~\eqref{\LocalName{eq:GradientMethod}} uses the following central
inequality (for $r_{k + 1} \DefinedEqual \Norm{x_{k + 1} - x_k}$ and
$d_k \DefinedEqual \Norm{x_k - x^*}$ with $x^*$ being a solution
of~\eqref{eq:MainProblem}):
\begin{multline*}
  f(x_k) + \InnerProduct{f'(x_k)}{x_{k + 1} - x_k} + \psi(x_{k + 1})
  +
  \frac{H_k}{2} r_{k + 1}^2
  +
  \frac{H_k}{2} d_{k + 1}^2
  \\
  \leq
  f(x_k) + \InnerProduct{f'(x_k)}{x^* - x_k} + \psi(x^*)
  +
  \frac{H_k}{2} d_k^2,
\end{multline*}
which is a simple consequence of the strong convexity of the objective function
in the auxiliary subproblem~\eqref{\LocalName{eq:GradientMethod}}
(c.f.\ \cref{th:OptimalityConditionForProximalPointStep}).
Rewriting now
$
  f(x_k) + \InnerProduct{f'(x_k)}{x_{k + 1} - x_k}
  =
  f(x_{k + 1}) - \BregmanDistance{f}(x_k, x_{k + 1})
$
using the Bregman distance, and estimating
$
  f(x_k) + \InnerProduct{f'(x_k)}{x^* - x_k} + \psi(x^*)
  \leq
  f(x^*) + \psi(x^*)
  =
  F^*
$
using the convexity of~$f$, we get
\begin{equation}
  \LocalLabel{eq:MainRecurrence}
  F(x_{k + 1}) - F^*
  +
  \frac{H_k}{2} d_{k + 1}^2
  \leq
  \frac{H_k}{2} d_k^2
  +
  \beta_{k + 1} - \frac{H_k}{2} r_{k + 1}^2,
\end{equation}
where
$
  \beta_{k + 1}
  \DefinedEqual
  f(x_{k + 1}) - f(x_k) - \InnerProduct{f'(x_k)}{x_{k + 1} - x_k}.
$

\subimport{./}{LineSearchApproach}
\subimport{./}{OurIdea}

%% file: GradientMethod/MainIdea/LineSearchApproach.tex
\subsubsection{Line-Search Approach}

The standard approach to proceed, pioneered by~\cite{nesterov2015universal},
is to choose the coefficient~$H_k$ in such a way that the error term
$\beta_{k + 1} - \frac{H_k}{2} r_{k + 1}^2$ in \cref{\LocalName{eq:MainRecurrence}}
is sufficiently small:
\begin{equation}
  \LocalLabel{eq:SmallErrorTerm}
  \Delta_k \DefinedEqual \beta_{k + 1} - \frac{H_k}{2} r_{k + 1}^2
  \leq
  \frac{\epsilon}{2}
\end{equation}
(for a certain fixed $\epsilon > 0$), and then divide both sides of
\cref{\LocalName{eq:MainRecurrence}} by~$H_k$ to get a telescopic recurrence:
\[
  \frac{1}{H_k} [F(x_{k + 1}) - F^*]
  +
  \frac{1}{2} d_{k + 1}^2
  \leq
  \frac{1}{2} d_k^2
  +
  \frac{\epsilon}{2 H_k}.
\]
Telescoping and dividing by
$S_k \DefinedEqual \sum_{i = 0}^{k - 1} \frac{1}{H_i}$, we get
\begin{equation}
  \LocalLabel{eq:ConvergenceGuaranteeForLineSearchMethod}
  F(x_k^*) - F^*
  \leq
  \frac{D_0^2}{2 S_k} + \frac{\epsilon}{2}
  \leq
  \frac{H_k^* D_0^2}{2 k} + \frac{\epsilon}{2},
\end{equation}
where $D_0 \DefinedEqual d_0 = \Norm{x_0 - x^*}$,
$H_k^* \DefinedEqual \max_{0 \leq i \leq k - 1} H_i$,
and $x_k^*$ is the ``best'' iterate:
\begin{equation}
  \label{eq:OutputOfGradientMethod}
  x_k^*
  \DefinedEqual
  \argmin \SetBuilder[\big]{f(x)}{x \in \Set{x_1, \ldots, x_k}}.
\end{equation}
(Alternatively, one could also define $x_k^*$ as the average of~$x_i$ with
weights~$\frac{1}{H_i}$.)
This gives us the convergence of the function residual to~$\epsilon$,
provided that $H_k^*$ is reasonably bounded from above (e.g., by a constant).

To ensure that \cref{\LocalName{eq:SmallErrorTerm}} is satisfied for a
sufficiently large $H_k$ and estimate the corresponding $H_k^*$, we start with
the observation that, by \cref{eq:UpperBoundOnBregmanDistanceViaHolderConstant},
$\beta_{k + 1} \leq \frac{L_{\nu}}{1 + \nu} r_{k + 1}^{1 + \nu}$ for any
$\nu \in \ClosedClosedInterval{0}{1}$, and hence
\begin{equation}
  \LocalLabel{eq:EstimatingBregmanErrorTerm}
  \Delta_k
  \leq
  \frac{L_{\nu}}{1 + \nu} r_{k + 1}^{1 + \nu} - \frac{H_k}{2} r_{k + 1}^2
  \leq
  \frac{
      (1 - \nu) L_{\nu}^{2 / (1 - \nu)}
  }{
      2 (1 + \nu) H_k^{(1 + \nu) / (1 - \nu)}
  }.
\end{equation}
(The final inequality follows by maximizing the expression in $r_{k + 1}$,
see \cref{th:UpperBoundOnRegularizedBregmanDistanceTerm};
for $\nu = 1$, the right-hand side should be understood
as $0$ if $H_k \geq L_{\nu}$ and $+\infty$ otherwise.)
Making the right-hand side of the above display $\leq \frac{\epsilon}{2}$, we see that
\cref{\LocalName{eq:SmallErrorTerm}} is satisfied
whenever $H_k \geq \bar{H}_{\nu}$, where
\[
  \bar{H}_{\nu}
  \DefinedEqual
  L_{\nu}^{2 / (1 + \nu)}
  \Bigl[ \frac{1 - \nu}{(1 + \nu) \epsilon} \Bigr]^{(1 - \nu) / (1 + \nu)}.
\]
Notice that $\nu = 1$ implies $\bar{H}_{\nu} \geq L_{\nu}$.
Since we do not know the best (= smallest) possible value of~$L_{\nu}$ over all
$\nu$, we cannot simply set
$
  H_k
  =
  \bar{H}_*
  \DefinedEqual
  \inf_{\nu \in \ClosedClosedInterval{0}{1}} \bar{H}_{\nu}
$.
However, we can let the line search estimate this value for us:
at each iteration, we start with a certain initial guess~$H_k'$ for~$H_k$
and then repeatedly double this value until
condition~\cref{\LocalName{eq:SmallErrorTerm}} is satisfied
(note that $x_{k + 1}$ depends on~$H_k$ and thus needs to be recomputed at
every iteration of the line search procedure).
Provided that the initial guess~$H_0'$ at the very first iteration is not
sufficiently large, e.g., $H_0' \leq \bar{H}_*$ and that $H_{k + 1}'$ is
chosen appropriately (e.g, $H_{k + 1}' = \frac{1}{2} H_k$), we then have the
guarantee that $H_k$ computed by the line search does not significantly
exceed the ``right'' value: $H_k \leq 2 \bar{H}_*$;
furthermore, the total number of line search iterations across all iterations
of the algorithm is reasonably bounded.
Substituting now this bound on $H_k$ into
\cref{\LocalName{eq:ConvergenceGuaranteeForLineSearchMethod}}
and looking at the number of iterations~$k$ that one needs to make
$\frac{H_k^* D_0^2}{k} \leq \epsilon$, we see that the outlined above
line-search method needs
\begin{equation}
  \LocalLabel{eq:ComplexityBoundForLineSearchMethod}
  \BigO\biggl(
    \inf_{\nu \in \ClosedClosedInterval{0}{1}}
    \frac{\bar{H}_{\nu} D_0^2}{\epsilon}
  \biggr)
  =
  \BigO\biggl(
    \inf_{\nu \in \ClosedClosedInterval{0}{1}}
    \biggl[ \frac{L_{\nu}}{\epsilon} \biggr]^{2 / (1 + \nu)} D_0^2
  \biggr)
\end{equation}
iterations to reach $F(x_k^*) - F^* \leq \epsilon$.

%% file: GradientMethod/MainIdea/OurIdea.tex
\subsubsection{Our Idea: How to Avoid Line Search}

The problem with the line-search approach, which makes it difficult to extend
the corresponding reasoning to the stochastic case, is that it creates a
dependency (correlation) between $x_k$ and $H_k$.
This does not allow us to use the unbiasedness of the stochastic gradient oracle
once we divide (the stochastic counterpart of)
\cref{\LocalName{eq:MainRecurrence}} by~$H_k$
(see \cref{sec:OutineOfAnalysisForStochasticGradientMethod}).

However, we can follow a different approach to convert
\cref{\LocalName{eq:MainRecurrence}} into a telescopic recurrence.
Specifically, let us replace the coefficient~$H_k$ in the left-hand side of
\cref{\LocalName{eq:MainRecurrence}} with~$H_{k + 1}$:
\[
  F(x_{k + 1}) - F^* + \frac{H_{k + 1}}{2} d_{k + 1}^2
  -
  \frac{H_k}{2} d_k^2
  \leq
  \beta_{k + 1} - \frac{H_k}{2} r_{k + 1}^2
  +
  \frac{1}{2} (H_{k + 1} - H_k) d_{k + 1}^2.
\]
As we can see, such an operation may introduce an additional error term
$\frac{1}{2} (H_{k + 1} - H_k) d_{k + 1}^2$ if we plan to increase our step-size
coefficient: $H_k \leq H_{k + 1}$ (which is a natural thing to do if it is
currently too small making the other error term
$\beta_{k + 1} - \frac{H_k}{2} r_{k + 1}^2$ too large).
Nevertheless, using \cref{as:BoundedFeasibleSet}, we can easily control this
additional error term and make it telescopic:
\begin{equation}
  \LocalLabel{eq:PreliminaryRecurrence}
  F(x_{k + 1}) - F^* + \frac{H_{k + 1}}{2} d_{k + 1}^2
  -
  \frac{H_k}{2} d_k^2
  \leq
  \beta_{k + 1} - \frac{H_k}{2} r_{k + 1}^2
  +
  \frac{1}{2} (H_{k + 1} - H_k) D^2.
\end{equation}

Our main idea now is to choose the next coefficient~$H_{k + 1}$ so that the two
error terms are balanced:
\begin{equation}
  \LocalLabel{eq:PreliminaryEquationForNextCoefficient}
  \frac{1}{2} (H_{k + 1} - H_k) D^2
  =
  \PositivePart[\Big]{\beta_{k + 1} - \frac{H_k}{2} r_{k + 1}^2},
\end{equation}
where we additionally put the positive part $\PositivePart{\cdot}$ to respect
the monotonicity relation~$H_k \leq H_{k + 1}$.
Recall that $\beta_{k + 1}$ and $r_{k + 1}$ depend only on $x_k$ and $x_{k + 1}$
(which themselves depend on $H_{k - 1}$ and $H_k$, see
\cref{\LocalName{eq:GradientMethod}}).
Thus, \cref{\LocalName{eq:PreliminaryEquationForNextCoefficient}} is a simple
linear equation for $H_{k + 1}$ which does not require any line search for
solving it.

Substituting our
choice~\eqref{\LocalName{eq:PreliminaryEquationForNextCoefficient}}
of~$H_{k + 1}$ into \cref{\LocalName{eq:PreliminaryRecurrence}}
(and using the fact that $\tau \leq \PositivePart{\tau}$ for
any~$\tau \in \RealField$), we arrive at the following simple telescopic
inequality:
\begin{equation}
  \LocalLabel{eq:FinalRecurrence}
  F(x_{k + 1}) - F^* + \frac{H_{k + 1}}{2} d_{k + 1}^2
  \leq
  \frac{H_k}{2} d_k^2 + (H_{k + 1} - H_k) D^2.
\end{equation}
Telescoping these inequalities, we get
\begin{equation}
  \LocalLabel{eq:PreliminaryEstimateForFunctionResidual}
  F(x_k^*) - F^*
  \leq
  \frac{1}{k} \Bigl[ \frac{H_0}{2} d_0^2 + (H_k - H_0) D^2 \Bigr]
  \leq
  \frac{H_k D^2}{k},
\end{equation}
where $x_k^*$ is the ``best'' iterate (see \cref{eq:OutputOfGradientMethod}).
(Alternatively, one could also define $x_k^* = \frac{1}{k} \sum_{i = 1}^k x_i$.)

The main question is how fast the coefficient~$H_k$ grows.
Following exactly the same argument as in
\cref{\LocalName{eq:EstimatingBregmanErrorTerm}}
(and using the fact that $\PositivePart{\cdot}$ is nondecreasing),
we can estimate the right-hand side of our balance
equation~\eqref{\LocalName{eq:PreliminaryEquationForNextCoefficient}}
and conclude that
$
  (H_{k + 1} - H_k) D^2
  \leq
  \frac{
    (1 - \nu) L_{\nu}^{2 / (1 - \nu)}
  }{
    (1 + \nu) H_k^{(1 + \nu) / (1 - \nu)}
  }.
$
(Assume, for simplicity, that $\nu < 1$;
to rigorously handle the case $\nu = 1$ we need a more careful argument.)
This is a certain recurrent inequality that we can use to estimate the rate of
growth of~$H_k$.
This would be especially simple if we had, say,
$2 H_{k + 1}^{(1 + \nu) / (1 - \nu)}$ instead of $H_k^{(1 + \nu) / (1 - \nu)}$:
\begin{equation}
  \LocalLabel{eq:RecurrentBoundForCoefficients}
  (H_{k + 1} - H_k) D^2
  \leq
  \frac{
    (1 - \nu) L_{\nu}^{2 / (1 - \nu)}
  }{
    2 (1 + \nu) H_{k + 1}^{(1 + \nu) / (1 - \nu)}
  }.
\end{equation}
Then, a simple integration argument
(see \cref{th:RecurrentInequalityWithPowerGrowth} with
$p = \frac{1 + \nu}{1 - \nu}$ for which $p + 1 = \frac{2}{1 - \nu}$)
would show that
\begin{equation}
  \LocalLabel{eq:FinalBoundForCoefficients}
  H_k \leq \frac{L_{\nu}}{D^{1 - \nu}} k^{(1 - \nu) / 2},
\end{equation}
provided that the initial step-size coefficient was chosen appropriately:
$H_0 = 0$.
(Note that this would not cause any problems for the
iteration~\eqref{\LocalName{eq:GradientMethod}} being well-defined since we
assume that $\EffectiveDomain \psi$ is a bounded set.)

Of course, we cannot argue that that our ``real'' version of
\cref{\LocalName{eq:RecurrentBoundForCoefficients}}
(the one with $H_k^{(1 + \nu) / (1 - \nu)}$ instead of
$2 H_{k + 1}^{(1 + \nu) / (1 - \nu)}$) implies the
desired \cref{\LocalName{eq:RecurrentBoundForCoefficients}}
(in fact, the relationship is exactly the opposite since $H_k \leq H_{k + 1}$).
However, we can slightly modify the reasoning we used to pass from
\cref{\LocalName{eq:PreliminaryRecurrence}} to
\cref{\LocalName{eq:PreliminaryEstimateForFunctionResidual}}
and the corresponding recurrent inequality for~$H_k$.
Specifically, we can rewrite the $-\frac{H_k}{2} r_{k + 1}^2$ term in the
right-hand side of \cref{\LocalName{eq:PreliminaryRecurrence}} as
$-\frac{H_{k + 1}}{2} r_{k + 1}^2 + \frac{1}{2} (H_{k + 1} - H_k) r_{k + 1}^2$,
and then upper bound $r_{k + 1} \leq D$.
As a result, we get
$\beta_{k + 1} - \frac{H_{k + 1}}{2} r_{k + 1}^2 + (H_{k + 1} - H_k) D^2$
in the right-hand of \cref{\LocalName{eq:PreliminaryRecurrence}},
and can now choose the coefficient~$H_{k + 1}$ using the following balance
equation instead of
\cref{\LocalName{eq:PreliminaryEquationForNextCoefficient}}:
\begin{equation}
  \LocalLabel{eq:SolvingEquationForRegularizationParameter}
  \boxed{
    (H_{k + 1} - H_k) D^2
    =
    \PositivePart[\Big]{\beta_{k + 1} - \frac{H_{k + 1}}{2} r_{k + 1}^2}.
  }
\end{equation}
Although this is no longer a linear equation in~$H_{k + 1}$, it always has a
unique solution $H_{k + 1} \geq H_k$, which can be easily computed:
if $\beta_{k + 1} \leq \frac{H_k}{2} r_{k + 1}^2$, then $H_{k + 1} = H_k$;
otherwise, $H_{k + 1}$ is the solution of the linear equation
$(H_{k + 1} - H_k) D^2 = \beta_{k + 1} - \frac{H_{k + 1}}{2} r_{k + 1}^2$
(see \cref{th:SolvingEquationForRegularizationParameter}).
Proceeding exactly is the same way as before, we get
\cref{\LocalName{eq:PreliminaryEstimateForFunctionResidual}} but with~$2 D^2$
instead of~$D^2$,
and, most importantly, the desired
\cref{\LocalName{eq:RecurrentBoundForCoefficients}}.
As a result, \cref{\LocalName{eq:FinalBoundForCoefficients}} indeed holds
and we get
\[
  F(x_k^*) - F^*
  \leq
  \frac{2 H_k D^2}{k}
  \leq
  \inf_{\nu \in \ClosedClosedInterval{0}{1}}
  \frac{2 L_{\nu} D^{1 + \nu}}{k^{(1 + \nu) / 2}},
\]
where the infimum is due to the fact that $\nu \in \ClosedClosedInterval{0}{1}$
was allowed to be arbitrary in our analysis.

%% file: GradientMethod/TheMethod.tex
\subsection{The Method}

Summarizing the outlined above considerations into a formal algorithmic scheme,
we arrive at \cref{alg:GradientMethod}.
This is essentially the classical (Composite) Gradient
Method~\eqref{\LocalName{eq:GradientMethod}} but equipped with our novel
step-size adjusting
rule~\eqref{\LocalName{eq:SolvingEquationForRegularizationParameter}}
(the formula for $H_{k + 1}$ at
\cref{GradientMethod::step:UpdateRegularizationParameter} is the explicitly
written solution of the balance
equation~\eqref{\LocalName{eq:SolvingEquationForRegularizationParameter}}).

\begin{algorithm}[tb]
  \caption{Universal Line-Search-Free Gradient Method}
  \label{alg:GradientMethod}
  \UsingNamespace{GradientMethod}
  \begin{algorithmic}[1]
    \State \textbf{Initialize}:
    $x_0 \in \EffectiveDomain \psi$, diameter $D > 0$, $H_0 \DefinedEqual 0$.

    \For{$k = 0, 1, \ldots$}
      \State
      \LocalLabel{step:ComputeSubgradient}
        Compute $g_k \in \Subdifferential f(x_k)$.

      \State
      \LocalLabel{step:ComputeNewPoint}
        $
          x_{k + 1}
          =
          \argmin_{x \in \EffectiveDomain \psi} \{
            \DualPairing{g_k}{x} + \psi(x) + \frac{H_k}{2} \Norm{x - x_k}^2
          \}
        $.

      \State
      \LocalLabel{step:UpdateRegularizationParameter}
        $
          H_{k + 1}
          \DefinedEqual
          H_k
          +
          \frac{
            \PositivePart{\beta_{k + 1} - \frac{1}{2} H_k r_{k + 1}^2}
          }{
            D^2 + \frac{1}{2} r_{k + 1}^2
          }
        $,
        where
        $r_{k + 1} \DefinedEqual \Norm{x_k - x_{k + 1}}$,
        $
          \beta_{k + 1}
          \DefinedEqual
          \BregmanDistanceWithSubgradient{f}{g_k}(x_k, x_{k + 1})
        $.
    \EndFor
  \end{algorithmic}
\end{algorithm}

\begin{restatable}{theorem}{thConvergenceRateOfGradientMethod}
  \label{th:ConvergenceRateOfGradientMethod}
  Let \cref{alg:GradientMethod} be applied to problem~\eqref{eq:MainProblem}
  under \cref{as:BoundedFeasibleSet,as:AtLeastOneHolderConstantIsFinite},
  and let $x_k^*$ be the ``best'' iterate as defined in
  \cref{eq:OutputOfGradientMethod}.
  Then, for all $k \geq 1$, we have
  \[
    F(x_k^*) - F^*
    \leq
    \inf_{\nu \in \ClosedClosedInterval{0}{1}}
    \frac{2 L_{\nu} D^{1 + \nu}}{k^{(1 + \nu) / 2}}.
  \]
  To reach $F(x_k^*) - F^* \leq \epsilon$ for any $\epsilon > 0$,
  it thus suffices to make
  $
    \inf_{\nu \in \ClosedClosedInterval{0}{1}}
    [\frac{2 L_{\nu}}{\epsilon}]^{2 / (1 + \nu)} D^2
  $
  iterations.
\end{restatable}

Comparing the efficiency bound from \cref{th:ConvergenceRateOfGradientMethod}
with the corresponding
bound~\eqref{\LocalName{eq:ComplexityBoundForLineSearchMethod}} for the
line-search method, we see that they are almost identical.
The only difference is that our method has the diameter of the feasible set $D$
instead of the initial distance~$D_0$.
However, as we show next, our method can be easily extended to stochastic
problems.

%% file: StochasticGradientMethod/Main.tex
\section{Universal Gradient Method for Stochastic Optimization}
\label{sec:GradientMethodForStochasticOptimization}

Now we assume that $f$ in problem~\eqref{eq:MainProblem} is accessible only via
the \emph{stochastic gradient oracle}~$\hat{g}$.
Formally, this is a pair $(g, \xi)$ consisting of a random variable~$\xi$ and a
mapping~$\Map{g}{\EffectiveDomain f \times \Image \xi}{\RealField^n}$
(with $\Image \xi$ being the image of~$\xi$).
When queried at a point~$x \in \EffectiveDomain \psi$, the oracle automatically
generates an independent copy~$\xi$ of its randomness, and then returns
$s = g(x, \xi)$ (notation: $s \SampledFrom \hat{g}(x)$)---a random estimate of
a subgradient of~$f$ at~$x$.

We make the following standard assumption on the oracle:

\begin{assumption}
  \label{as:StochasticGradientOracle}
  The function~$f$ in problem~\eqref{eq:MainProblem} is accessible only via
  an unbiased stochastic gradient oracle~$\hat{g} = (g, \xi)$ with bounded
  variance:
  \begin{gather}
    \label{eq:StochasticGradientIsUnbiased}
    f'(x)
    \DefinedEqual
    \Expectation_{\xi}[g(x, \xi)]
    \in
    \Subdifferential f(x),
    \\
    \label{eq:VarianceOfStochasticGradient}
    \sigma^2
    \DefinedEqual
    \sup_{x \in \EffectiveDomain \psi}
    \Expectation_{\xi}[\DualNorm{g(x, \xi) - f'(x)}^2]
    <
    +\infty.
  \end{gather}
\end{assumption}

\subimport{./}{TheMethod}
\subimport{./}{MainIdea}
\subimport{./}{ComparisonWithAdagrad}

%% file: StochasticGradientMethod/MainIdea.tex
\subsection{Main Idea and Outline of Analysis}
\label{sec:OutineOfAnalysisForStochasticGradientMethod}

\UsingNamespace{OutlineOfAnalysisForStochasticGradientMethod}

Let us briefly explain the motivation behind the specific formula for
$\hat{\beta}_{k + 1}$ in \cref{alg:StochasticGradientMethod} and sketch the
corresponding convergence analysis.
The formal proof with all the details can be found in
\cref{sec:ProofForStochasticGradientMethod}.

From the definition of~$x_{k + 1}$, it follows that
(\cref{th:OptimalityConditionForProximalPointStep})
\begin{multline*}
  f(x_k) + \InnerProduct{g_k}{x_{k + 1} - x_k} + \psi(x_{k + 1})
  +
  \frac{H_k}{2} r_{k + 1}^2 + \frac{H_k}{2} d_{k + 1}^2
  \\
  \leq
  f(x_k) + \InnerProduct{g_k}{x^* - x_k} + \psi(x^*) + \frac{H_k}{2} d_k^2,
\end{multline*}
where $r_{k + 1} \DefinedEqual \Norm{x_{k + 1} - x_k}$
and $d_k \DefinedEqual \Norm{x_k - x^*}$.

Observe that
$
  \Expectation_{\xi_k} [f(x_k) + \InnerProduct{g_k}{x^* - x_k}]
  =
  f(x_k) + \InnerProduct{f'(x_k)}{x^* - x_k}
  \leq
  f(x^*)
$,
where $\xi_k$ is the oracle's randomness defining~$g_k \equiv g(x_k, \xi_k)$.
However, if we attempted to follow the line-search idea from
\cref{sec:MainIdeaForGradientMethod} by first dividing both sides in the
previous display by~$H_k$, then we would not be able to use the oracle's
unbiasedness as $H_k$ and~$x_k$ would depend on each other.

Nevertheless, our line-search-free idea still works.
Specifically, passing to expectations in the above display and using the lower
bound on $f(x^*)$ from the previous paragraph, and then rearranging, we obtain
\begin{equation}
  \LocalLabel{eq:PreliminaryRecurrence}
  \Expectation \Bigl[
    F(x_{k + 1}) - F^*
    +
    \frac{H_{k + 1}}{2} d_{k + 1}^2 - \frac{H_k}{2} d_k^2
  \Bigr]
  \leq
  \Expectation \Bigl[
    \beta_{k + 1} - \frac{H_{k + 1}}{2} r_{k + 1}^2
    +
    (H_{k + 1} - H_k) D^2
  \Bigr],
\end{equation}
where
$
  \beta_{k + 1}
  \DefinedEqual
  f(x_{k + 1}) - f(x_k) - \InnerProduct{g_k}{x_{k + 1} - x_k}
$,
and the $(H_{k + 1} - H_k) D^2$ term corresponds to the upper bound on
$\frac{1}{2} (H_{k + 1} - H_k) (d_{k + 1}^2 + r_{k + 1}^2)$.

The problem is that we cannot compute $\beta_{k + 1}$ since it involves the
exact function values~$f(x_{k + 1})$ and~$f(x_k)$.
However, we may replace it with an appropriate stochastic approximation.
Indeed, for our goals it suffices to know some $\hat{\beta}_{k + 1}$ which is an
upper estimate for~$\beta_{k + 1}$ in expectation:
$\Expectation \beta_{k + 1} \leq \Expectation \hat{\beta}_{k + 1}$.
To get an appropriate $\hat{\beta}_{k + 1}$, we could, in principle, ask the
oracle to provide not only stochastic gradients but also stochastic function
values.
However, this would require imposing extra requirements for the oracle.

Instead, we take another, simpler, approach.
By the convexity of~$f$, we can estimate
\begin{equation}
  \LocalLabel{eq:EstimateForBregmanDistance}
  \beta_{k + 1}
  \leq
  \InnerProduct{f'(x_{k + 1}) - g_k}{x_{k + 1} - x_k}
  =
  \Expectation_{\xi_{k + 1}}[\hat{\beta}_{k + 1}],
\end{equation}
where
$
  \hat{\beta}_{k + 1}
  \DefinedEqual
  \InnerProduct{g_{k + 1} - g_k}{x_{k + 1} - x_k}
$
can be calculated in the algorithm and $\xi_{k + 1}$ is the oracle's randomness
defining~$g_{k + 1} \equiv g(x_{k + 1}, \xi_{k + 1})$.
It is important for the final identity that $\xi_{k + 1}$ is generated after
$x_k$ and $x_{k + 1}$.

This leads us to the balance equation
\begin{equation}
  \LocalLabel{eq:BalanceEquation}
  \boxed{
    (H_{k + 1} - H_k) D^2
    =
    \PositivePart{\hat{\beta}_{k + 1} - \tfrac{1}{2} H_{k + 1} r_{k + 1}^2},
  }
\end{equation}
whose solution is given at
\cref{StochaticGradientMethod::step:UpdateRegularizationParameter}
in \cref{alg:StochasticGradientMethod}.

Taking into account our balance equation and
\cref{\LocalName{eq:EstimateForBregmanDistance}}, we obtain exactly the same
simple-to-telescope inequality as
\cref{MainIdeaForGradientMethod::eq:FinalRecurrence} (valid in expectation)
which then leads to
\begin{equation}
  \LocalLabel{eq:PreliminaryConvergenceRateBound}
  \Expectation[F(\bar{x}_k)] - F^* \leq \frac{2 \Expectation[H_k] D^2}{k}.
\end{equation}
The rest of the analysis focuses on estimating the (expected) rate of growth
of~$H_k$.
The key idea is that we can estimate
\[
  \hat{\beta}_{k + 1}
  =
  \InnerProduct{f'(x_{k + 1}) - f'(x_k) + \Delta_{k + 1}}{x_{k + 1} - x_k}
  \leq
  L_{\nu} r_{k + 1}^{1 + \nu} + \sigma_{k + 1} r_{k + 1},
\]
where
$f'(x_k) \DefinedEqual \Expectation_{\xi_k}[g_k] \in \Subdifferential f(x_k)$
and $\Delta_{k + 1} \DefinedEqual \delta_{k + 1} - \delta_k$ with
$\delta_k \DefinedEqual g_k - f'(x_k)$ being the error of the stochastic
gradient (such that $\Expectation \Norm{\delta_k}^2 \leq \sigma^2$),
and $\sigma_{k + 1} \DefinedEqual \Norm{\Delta_{k + 1}}$.
This gives
\[
  (H_{k + 1} - H_k) D^2
  =
  \PositivePart{
    L_{\nu} r_{k + 1}^{1 + \nu}
    +
    \sigma_{k + 1} r_{k + 1}
    -
    \tfrac{1}{2} H_{k + 1} r_{k + 1}^2
  }.
\]
Eliminating $r_{k + 1}$ from this inequality
(by maximizing the right-hand side in this variable),
we get a certain recurrence for~$H_{k + 1}$, which is similar to
\cref{MainIdeaForGradientMethod::eq:RecurrentBoundForCoefficients} but with an
additional $\frac{\sigma_{k + 1}^2}{H_{k + 1}}$ term in the right-hand side;
carefully analyzing the resulting recurrence
(see \cref{th:EstimatingGrowthRateOfRegularizationParameterInStochasticCase}),
we get
$
  H_k
  \leq
  \BigO\bigl(
    \frac{L_{\nu}}{D^{1 - \nu}} k^{(1 - \nu) / 2}
    +
    \frac{1}{D} ( \sum_{i = 1}^k \sigma_i^2 )^{1 / 2}
  \bigr).
$
This gives us
\[
  \Expectation[H_k]
  \leq
  \BigO\Bigl(
    \frac{L_{\nu}}{D^{1 - \nu}} k^{(1 - \nu) / 2} + \frac{\sigma}{D} \sqrt{k}
  \Bigr)
\]
after taking expectations.
Substituting this bound into
\cref{\LocalName{eq:PreliminaryConvergenceRateBound}}, we get the convergence
result from \cref{th:ConvergenceRateOfStochasticGradientMethod}.

%% file: StochasticGradientMethod/ComparisonWithAdagrad.tex
\subsection{Comparison with AdaGrad-type Methods}

Let us compare the step-size coefficient~$H_k$ from
\cref{alg:StochasticGradientMethod} with that of AdaGrad-type methods.
Denote
\[
  \gamma_{k + 1} \DefinedEqual \DualNorm{g_{k + 1} - g_k}.
\]
From the definitions of~$\hat{\beta}_{k + 1}$ and~$r_{k + 1}$, it follows that
$
  \hat{\beta}_{k + 1} - \frac{H_{k + 1}}{2} r_{k + 1}^2
  \leq
  \gamma_{k + 1} r_{k + 1} - \frac{H_{k + 1}}{2} r_{k + 1}^2
  \leq
  \frac{\gamma_{k + 1}^2}{2 H_{k + 1}}
$.
Substituting this into the balance
equation~\eqref{\LocalName{eq:BalanceEquation}}
(using the monotonicity of~$\PositivePart{\cdot}$),
we get
\begin{equation}
  \LocalLabel{eq:InequalityForStepSizeCoefficientViaGradientNorm}
  (H_{k + 1} - H_k) D^2 \leq \frac{\gamma_{k + 1}^2}{2 H_{k + 1}}.
\end{equation}
From this and $H_0 = 0$, it follows that
(see \cref{th:RecurrentInequalityWithPowerGrowth})
\begin{equation}
  \LocalLabel{eq:RelationToAdagradStepSizeCoefficient}
  H_k
  \leq
  H_k'
  \DefinedEqual
  \frac{1}{D} \Bigl( \sum_{i = 1}^k \gamma_i^2 \Bigr)^{1 / 2}.
\end{equation}
Note that $H_k'$ is the step-size coefficient used by a variety of other
AdaGrad-type algorithms\footnote{
  The classical AdaGrad uses $\gamma_i = \DualNorm{g_i}$ but such a choice does
  not work properly for smooth constrained optimization when
  $\nabla f(x^*) \neq 0$.
}
(see, e.g.,~\cite{kavis2019unixgrad,ene2022adaptive}).
Thus, the ``step size'' $\frac{1}{H_k}$ in our algorithm is at least as large as
$\frac{1}{H_k'}$ used by AdaGrad.

In fact, the theoretical reasoning we used in
\cref{sec:OutineOfAnalysisForStochasticGradientMethod} to arrive at our formula
for the step-size coefficient, can be seen as a more precise theoretical
analysis of the Stochastic Gradient Method with adaptive step sizes.
Specifically, coming back to our preliminary
recurrence~\eqref{\LocalName{eq:PreliminaryRecurrence}}, we see that AdaGrad
first estimates
$
  \hat{\beta}_{k + 1} - \frac{H_{k + 1}}{2} r_{k + 1}^2
  \leq
  \frac{\gamma_{k + 1}^2}{2 H_{k + 1}}
$
and only then attempts to balance the terms.
This corresponds to the idea of choosing the coefficient~$H_{k + 1}$ in such a
way that \cref{\LocalName{eq:InequalityForStepSizeCoefficientViaGradientNorm}}
becomes an identity
(and then we not only have
\cref{\LocalName{eq:RelationToAdagradStepSizeCoefficient}}
but also the similar lower bound $H_k \geq \frac{1}{\sqrt{2}} H_k'$
(see \cref{th:RecurrentInequalityWithQuadraticGrowth}),
which means that $H_k = \Theta(H_k')$).
In contrast, our reasoning suggests that this extra estimation step is
unnecessary.

%% file: StochasticFastGradientMethod.tex
\section{Universal Fast Gradient Method for Stochastic Optimization}
\label{sec:FastGradientMethodForStochasticOptimization}

\begin{algorithm}[tb]
  \caption{Universal Stochastic Fast Gradient Method}
  \label{alg:StochasticFastGradientMethod}
  \UsingNamespace{StochasticFastGradientMethod}
  \begin{algorithmic}[1]
    \State \textbf{Initialize:}
    $x_0 = v_0 \in \EffectiveDomain \psi$, $D > 0$,
    $H_0 \DefinedEqual A_0 \DefinedEqual 0$.

    \For{$k = 0, 1, \ldots$}
      \State
      \LocalLabel{step:ChooseScalingCoefficient}
      $a_{k + 1} \DefinedEqual k + 1$, \quad
      $A_{k + 1} \DefinedEqual A_k + a_{k + 1} \ (> 0)$.

      \State
      \LocalLabel{step:ComputeIntermediatePoint}
      $
        y_k
        \DefinedEqual
        \frac{A_k}{A_{k + 1}} x_k + \frac{a_{k + 1}}{A_{k + 1}} v_k
      $, \quad
      $g_k^y \SampledFrom \hat{g}(y_k)$.

      \State
      \LocalLabel{step:UpdateProxCenter}
      $
        v_{k + 1}
        =
        \argmin_{x \in \EffectiveDomain \psi} \{
          a_{k + 1} [\DualPairing{g_k^y}{x} + \psi(x)]
          +
          \frac{H_k}{2} \Norm{x - v_k}^2
        \}.
      $

      \State
      \LocalLabel{step:ComputeNewPoint}
      $
        x_{k + 1}
        \DefinedEqual
        \frac{A_k}{A_{k + 1}} x_k + \frac{a_{k + 1}}{A_{k + 1}} v_{k + 1}
      $.

      \State
      \LocalLabel{step:UpdateRegularizationParameter}
      $
        H_{k + 1}
        \DefinedEqual
        H_k
        +
        \frac{
          \PositivePart{
            A_{k + 1} \hat{\beta}_{k + 1} - \frac{1}{2} H_k r_{k + 1}^2
          }
        }{
          D^2 + \frac{1}{2} r_{k + 1}^2
        }
      $,\par
      where
      $r_{k + 1} = \Norm{v_{k + 1} - v_k}$,
      $
        \hat{\beta}_{k + 1}
        =
        \DualPairing{g_{k + 1}^x - g_k^y}{x_{k + 1} - y_k}
      $
      with $g_{k + 1}^x \SampledFrom \hat{g}(x_{k + 1})$.
    \EndFor
  \end{algorithmic}
\end{algorithm}

We now present, in \cref{alg:StochasticFastGradientMethod}, the accelerated
version of our universal stochastic method for solving
problem~\eqref{eq:MainProblem}.
This algorithm is essentially one of the standard variants of the Fast Gradient
Method known as the Method of Similar Triangles (see, e.g., Section~6.1.3
in~\cite{nesterov2018lectures}), which uses stochastic gradients instead of the
exact ones and is equipped with our novel rule for adjusting the step-size
coefficient~$H_k$.
The algorithm enjoys the following efficiency estimate
(see \cref{sec:ProofForStochasticFastGradientMethod} for the proof):

\begin{restatable}{theorem}{thConvergenceRateOfStochasticFastGradientMethod}
  \label{th:ConvergenceRateOfStochasticFastGradientMethod}
  Let \cref{alg:StochasticFastGradientMethod} be applied to
  problem~\eqref{eq:MainProblem} under
  \cref{%
    as:BoundedFeasibleSet,%
    as:AtLeastOneHolderConstantIsFinite,%
    as:StochasticGradientOracle%
  }.
  Then, for all $k \geq 1$,
  \[
    \Expectation[F(x_k)] - F^*
    \leq
    \inf_{\nu \in \ClosedClosedInterval{0}{1}}
    \frac{32 L_{\nu} D^{1 + \nu}}{k^{(1 + 3 \nu) / 2}}
    +
    \frac{8 \sigma D}{\sqrt{3 k}}.
  \]
  Consequently, it suffices to make
  $
    \BigO\bigl(
      \inf_{\nu \in \ClosedClosedInterval{0}{1}}
      [\frac{L_{\nu} D^{1 + \nu}}{\epsilon}]^{2 / (1 + 3 \nu)}
      +
      \frac{\sigma^2 D^2}{\epsilon^2}
    \bigr)
  $
  oracle calls to reach $\Expectation[F(x_k)] - F^* \leq \epsilon$
  for any $\epsilon > 0$.
\end{restatable}

In the deterministic case ($\sigma = 0$), the efficiency bound from
\cref{th:ConvergenceRateOfStochasticFastGradientMethod} coincides with that of
the Universal Fast Gradient Method from~\cite{nesterov2015universal}.
It is worth mentioning that, in this case, instead of the symmetrized Bregman
distance, we can use the standard one,
$
  \hat{\beta}_{k + 1}
  =
  \BregmanDistanceWithSubgradient{f}{g_k^y}(y_k, x_{k + 1})
$,
in
\cref{alg:StochasticFastGradientMethod}, and get similar convergence estimates
but with slightly better absolute constants
(see \cref{th:ConvergenceRateOfFastGradientMethod}).

%% file: Experiments/Main.tex
\section{Experiments}

In this section, we present some preliminary computational experiments for the
proposed methods.

\subimport{./}{ConvexOptimization}
\subimport{./}{NeuralNetworks}

%% file: Experiments/ConvexOptimization.tex
\subsection{Convex Optimization}
\label{sec:exp_convex}

\begin{figure}[tb]
  \begin{center}
    \includegraphics[width=0.8\linewidth]{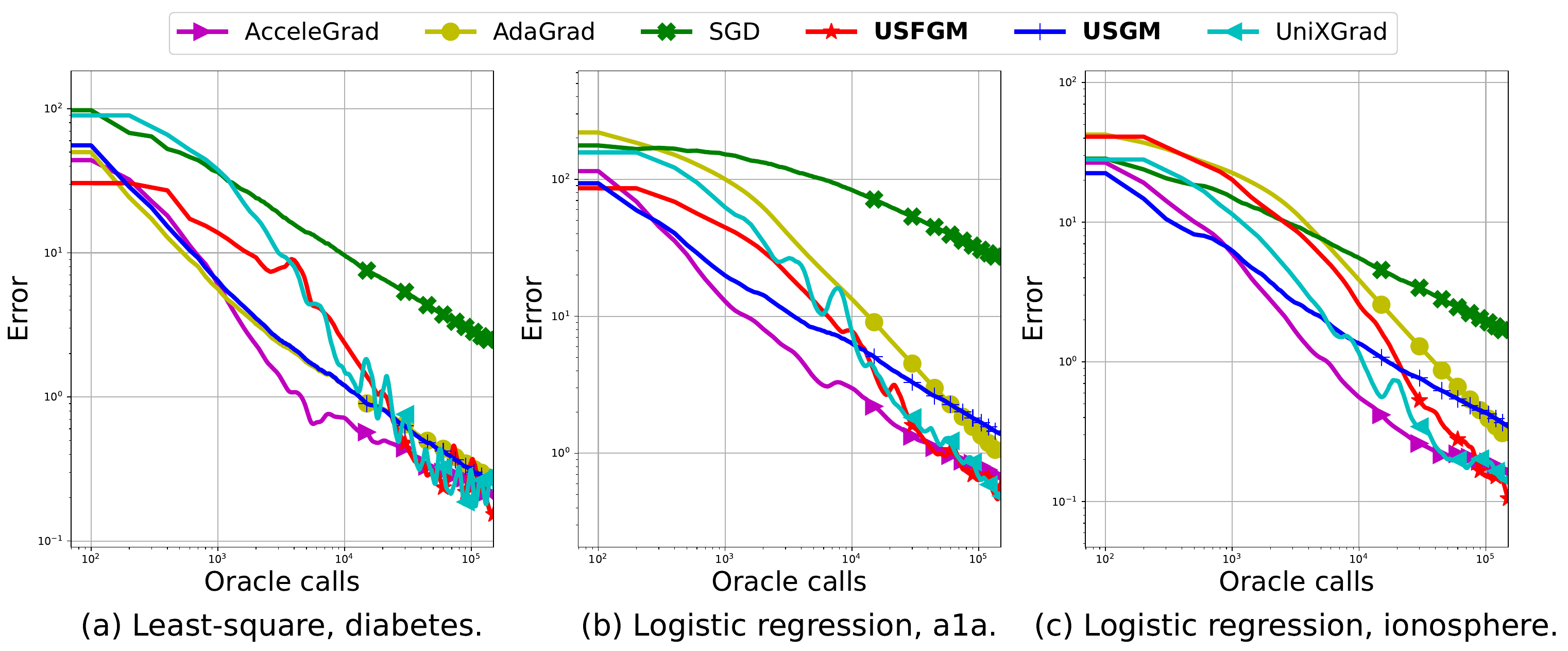}
  \end{center}
  \caption{Comparison of different stochastic algorithms on convex problems.}
  \label{fig:exp1}
\end{figure}

\paragraph{Least Squares.}

Let us consider the following problem:
\begin{equation}
  \label{equ:ls}
  \min_{x \in \RealField^n} \SetBuilder[\Big]{
    F(x) \DefinedEqual \frac{1}{2} \Norm{A x - b}^2
  }{\Norm{x} \leq 1},
\end{equation}
where $A \in \RealField^{m \times n}$, $b \in \RealField^m$.
We run the experiment on real-world diabetes dataset from LIBSVM\footnote{%
  \url{https://www.csie.ntu.edu.tw/~cjlin/libsvmtools/datasets/}%
}.
In the stochastic setting, we run the proposed
USGM (\cref{alg:StochasticGradientMethod}) and
USFGM (\cref{alg:StochasticFastGradientMethod}),
and compare them against SGD, AdaGrad, UnixGrad~\parencite{kavis2019unixgrad},
and AcceleGrad~\parencite{levy2018online}.
The result in \cref{fig:exp1}.(a)
shows that the proposed method attains a convergence rate comparable to AdaGrad
and AcceleGrad.

\paragraph{Logistic regression.}

We also focus on the logistic loss:
\[
  \min_{x \in \RealField^n} \SetBuilder[\Big]{
    F(x)
    \DefinedEqual
    \sum_{i = 1}^m \log\bigl( 1 + \exp(-b_i \InnerProduct{a_i}{x}) \bigr)
  }{\Norm{x} \leq 1},
\]
where $a_i \in \RealField^n$ is the feature vector and $b_i \in \Set{0, 1}$
is the label.
We run the experiment on the a1a and ionosphere datasets from LIBSVM.
The remaining setup is the same as in the case of Least-Squares.
We present the result in \cref{fig:exp1}.(b-c), where ASUGM and SUGM are
slightly faster than AdaGrad while performing similarly to UniXGrad.

%% file: Experiments/NeuralNetworks.tex
\subsection{Non-Convex Neural Networks Training}
\label{sec:exp_nn}

\begin{figure}[tb]
  \begin{center}
    \includegraphics[width=0.6\linewidth]{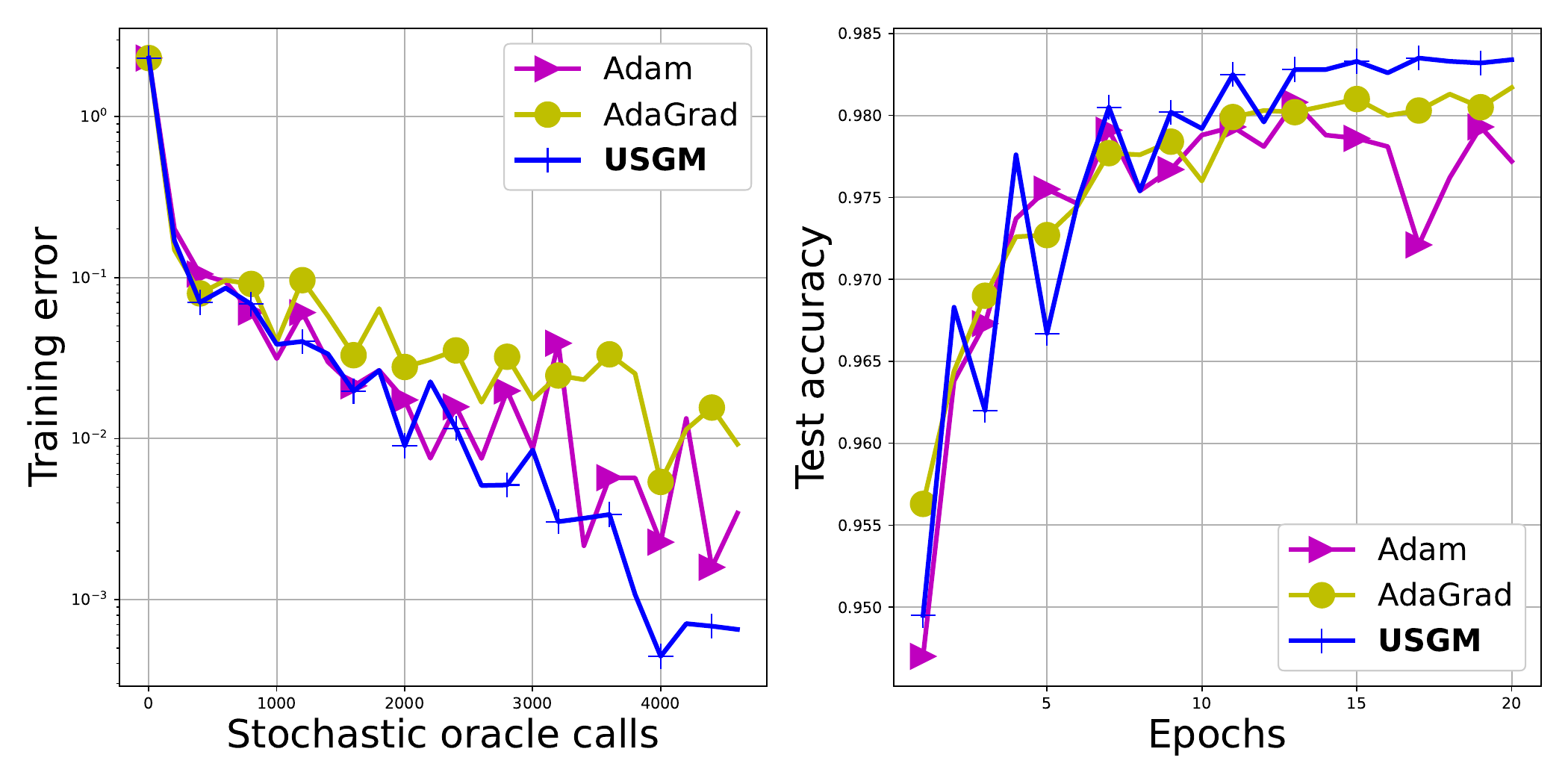}
  \end{center}
  \caption{%
    Comparison between the proposed universal stochastic gradient method, Adam,
    and AdaGrad in neural network training.%
  }
  \label{fig:exp2}
\end{figure}

We now show that the proposed method can also be applied to non-convex neural
network training.
Specifically, we focus on classification tasks with the cross-entropy loss on
MNIST dataset.
A three-layer fully connected networks with layer
dimensions~$[28 \times 28, 256, 256, 10]$ and ReLU activation function are
selected.
We compare the proposed method against AdaGrad and Adam.
We select the mini-batch size as~$256$.
Step-size of each method is tuned by a parameter sweep over
$\Set{10, 1, 0.1, 0.01, 0.001, 0.0001}$.
The diameter of the proposed method is tuned by sweeping over
$\Set{50, 35, 20, 10, 5}$.
We present the result in \cref{fig:exp2}, where we can see the proposed
stochastic universal gradient method can solve non-convex problems as well.

%% file: Conclusion.tex
\section{Conclusion}

We have proposed first universal gradient methods that are provably adaptive
simultaneously to the noise level in the gradient feedback, the Hölder exponent
and the associated Hölder constant of the objective function, while achieving
optimal efficiency bounds.
Unlike the majority of the works on adaptive methods relying on
AdaGrad step-size constructions, our algorithm design is inspired by the
line-search approach.
We have proposed a nonlinear balance equation for updating the step-size
coefficient, which results in a tighter analysis of adaptive stochastic
gradient algorithms compared to existing AdaGrad methods.

Note that our analysis exploits the fact that the feasible set has the
bounded diameter~$D$, knowledge of which is available to our algorithms.
While this assumption may seem rather restrictive, it is nevertheless quite
similar to the classical assumption on the knowledge of an upper bound~$R_0$
on the distance~$\Norm{x_0 - x^*}$ from the initial point to the minimizer.
Indeed, if we know $R_0$, we can easily convert our original
problem~\eqref{eq:MainProblem} into an equivalent one,
$\min_{x \in \EffectiveDomain \psi_D} [f(x) + \psi_D(x)]$,
where $\psi_D$ is the restriction of $\psi$ onto the ball
$
  B_0
  \DefinedEqual
  \SetBuilder{x \in \RealField^n}{\Norm{x - x_0} \leq R_0}
$
or, in other words, the sum of~$\psi$ and the indicator function of~$B_0$.
For this new problem, we can run our methods with diameter $D = 2 R_0$.
The only detail that one needs to address is how to compute the proximal-point
step for the function $\psi_D$ via the corresponding operation for $\psi$.
But this is usually not difficult and requires solving a certain one-dimensional
nonlinear equation, which can be done very efficiently by Newton's method
(at no extra queries to the stochastic gradient oracle).
In some special cases, this equation can even be solved analytically, e.g., when
the original problem is unconstrained, one simply needs to perform the
projection on the $B_0$ ball.
Nonetheless, it is still an interesting open question whether the same type
of results can be obtained without this (somewhat artificial) replacement of the
original problem.
More importantly, it would be interesting to obtain parameter-free versions
of our algorithms similar to~\cite{carmon2022making,ivgi2023dog}, which could
work with a sufficiently loose approximation of~$R_0$.

Having that said, we would like to make a few remarks regarding the technical
challenges involved in the design of optimal universal methods for stochastic
optimization.
Note that the existing accelerated adaptive methods for minimizing smooth
convex functions
\parencite{levy2018online,kavis2019unixgrad,joulani2020simpler}, which assume
bounded domains and use AdaGrad-inspired step-sizes, do not trivially extend to
the Hölder class of functions.
Essentially, they rely on the knowledge that the objective function is either
Lipschitz smooth ($\nu = 1$) or Lipschitz continuous ($\nu = 0$), and the
analysis is not directly compatible with intermediate values of
$\nu \in \OpenOpenInterval{0}{1}$.
Our approach is based on different, more suitable, techniques and yields a new
adaptive step-size schedule that enables fast universal rates in the
noisy setting.

%% file: Appendix/Main.tex
\subimport{./}{ProofForGradientMethod}
\subimport{./}{ProofForStochasticGradientMethod}
\subimport{./}{ProofForStochasticFastGradientMethod}
\subimport{./}{FastGradientMethod}
\subimport{./}{AuxiliaryResults}

\subimport{./}{AdditionalRelatedWork}
\subimport{AdditionalExperiments/}{Main}

%% file: Appendix/ProofForGradientMethod.tex
\section{Proof of Theorem~\ref{th:ConvergenceRateOfGradientMethod}}

\thConvergenceRateOfGradientMethod*

\begin{proof}
  \UsingNestedNamespace{ConvergenceRateOfGradientMethod}{Proof}

  \ProofPart

  We are going to prove that, for any $k \geq 1$,
  \begin{equation}
    \UsingNamespace{ConvergenceRateOfGradientMethod}
    \LocalLabel{eq:Result}
    F(x_k^*) - F^*
    \leq
    F(x_k^*) - \Phi_k^*
    \leq
    \frac{1}{k} \sum_{i = 1}^k F(x_i) - \Phi_k^*
    \leq
    \frac{2 H_k D^2}{k}
    \leq
    \inf_{\nu \in \ClosedClosedInterval{0}{1}}
    \frac{2 L_{\nu} D^{1 + \nu}}{k^{(1 + \nu) / 2}}.
  \end{equation}
  where
  \begin{equation}
    \label{eq:EstimatingFunctionInGradientMethod}
    \Phi_k^*
    \DefinedEqual
    \min_{x \in \EffectiveDomain \psi} \Bigl\{
      \Phi_k(x)
      \DefinedEqual
      \frac{1}{k} \sum_{i = 0}^{k - 1} [f(x_i) + \DualPairing{g_i}{x - x_i}]
      +
      \psi(x)
    \Bigr\}
    \quad (\leq F^*).
  \end{equation}
  (The inequality follows from the fact that $g_i \in \Subdifferential f(x_i)$
  for all $i \geq 0$ and \cref{eq:MainProblem}.)

  Note that the function $\Phi_k$, defined in
  \cref{eq:EstimatingFunctionInGradientMethod}, is the sum of
  an affine function and~$\psi$.
  Since $\psi$ is simple by our assumptions, we can easily compute its minimal
  value~$\Phi_k^*$.
  This value allows us to compute the quantities
  $\epsilon_k^* \DefinedEqual F(x_k^*) - \Phi_k^*$ and
  $\bar{\epsilon}_k \DefinedEqual \frac{1}{k} \sum_{i = 1}^k F(x_i) - \Phi_k^*$,
  appearing in \cref{ConvergenceRateOfGradientMethod::eq:Result},
  and thus equip \cref{alg:GradientMethod} with a reliable stopping criterion
  $\epsilon_k^* \leq \epsilon$ (or $\bar{\epsilon}_k \leq \epsilon$) which guarantees that
  $F(x_k^*) - F^* \leq \epsilon$ for some given $\epsilon > 0$.

  The first inequality in \cref{\MasterName{eq:Result}} follows from
  \cref{eq:EstimatingFunctionInGradientMethod}.
  The second one follows from the definition of~$x_k^*$ in
  \cref{eq:OutputOfGradientMethod}.

  \ProofPart

  Let us prove the third inequality in \cref{\MasterName{eq:Result}}.

  For each $k \geq 0$, let $\Map{\zeta_k}{\RealField^n}{\RealFieldPlusInfty}$
  be the function
  \begin{equation}
    \LocalLabel{eq:PartOfEstimatingFunction}
    \zeta_k(x)
    \DefinedEqual
    f(x_k) + \DualPairing{g_k}{x - x_k} + \psi(x).
  \end{equation}

  Let $k \geq 0$ and $x \in \EffectiveDomain \psi$ be arbitrary.
  By the definition of $x_{k + 1}$ at
  \cref{GradientMethod::step:ComputeNewPoint} and
  \cref{th:OptimalityConditionForProximalPointStep}, we have
  \begin{equation}
    \LocalLabel{eq:OptimalityCondition}
    \zeta_k(x) + \tfrac{1}{2} H_k \Norm{x - x_k}^2
    \geq
    \zeta_k(x_{k + 1}) + \tfrac{1}{2} H_k \Norm{x_{k + 1} - x_k}^2
    +
    \tfrac{1}{2} H_k \Norm{x - x_k}^2.
  \end{equation}
  According to
  \cref{%
    \LocalName{eq:PartOfEstimatingFunction},%
    eq:BregmanDistance,eq:MainProblem%
  },
  \[
    \zeta_k(x_{k + 1})
    =
    f(x_k) + \DualPairing{g_k}{x_{k + 1} - x_k} + \psi(x_{k + 1})
    =
    F(x_{k + 1}) - \BregmanDistanceWithSubgradient{f}{g_k}(x_k, x_{k + 1}).
  \]
  Substituting this into \cref{\LocalName{eq:OptimalityCondition}} and taking
  into account the definitions of~$r_{k + 1}$ and~$\beta_{k + 1}$, we get
  \begin{align*}
    \zeta_k(x) + \tfrac{1}{2} H_k \Norm{x - x_k}^2
    &\geq
    F(x_{k + 1})
    -
    \BregmanDistanceWithSubgradient{f}{g_k}(x_k, x_{k + 1})
    +
    \tfrac{1}{2} H_k \Norm{x_{k + 1} - x_k}^2
    +
    \tfrac{1}{2} H_k \Norm{x - x_k}^2
    \\
    &=
    F(x_{k + 1})
    -
    \beta_{k + 1} + \tfrac{1}{2} H_k r_{k + 1}^2
    +
    \tfrac{1}{2} H_k \Norm{x - x_k}^2.
  \end{align*}
  Consequently,
  \begin{equation}
    \LocalLabel{eq:PreliminaryRecursion}
    \begin{aligned}
      F(x_{k + 1}) + \tfrac{1}{2} H_{k + 1} \Norm{x - x_{k + 1}}^2
      &\leq
      \zeta_k(x) + \tfrac{1}{2} H_k \Norm{x - x_k}^2
      +
      [\beta_{k + 1} - \tfrac{1}{2} H_{k + 1} r_{k + 1}^2]
      \\
      &\qquad
      +
      \tfrac{1}{2} (H_{k + 1} - H_k) (\Norm{x - x_{k + 1}}^2 + r_{k + 1}^2).
    \end{aligned}
  \end{equation}

  Note that, by construction, $H_k \leq H_{k + 1}$.
  Also, in view of \cref{as:BoundedFeasibleSet} (and the fact that $x_i \in
  \EffectiveDomain \psi$ for any $i \geq 0$), we have $r_{k + 1} \leq D$ and
  $\Norm{x - x_{k + 1}} \leq D$.
  Therefore,
  \begin{equation}
    \LocalLabel{eq:UpperBoundOnDistanceTerm}
    \tfrac{1}{2} (H_{k + 1} - H_k) (\Norm{x - x_{k + 1}}^2 + r_{k + 1}^2)
    \leq
    (H_{k + 1} - H_k) D^2.
  \end{equation}
  At the same time, by the definition of $H_{k + 1}$ at
  \cref{GradientMethod::step:UpdateRegularizationParameter} and
  \cref{th:SolvingEquationForRegularizationParameter}, it satisfies the
  following equation:
  \begin{equation}
    \LocalLabel{eq:EquationForRegularizationParameter}
    (H_{k + 1} - H_k) D^2
    =
    \PositivePart{\beta_{k + 1} - \tfrac{1}{2} H_{k + 1} r_{k + 1}^2}.
  \end{equation}
  Therefore,
  \begin{equation}
    \LocalLabel{eq:UpperBoundOnBregmanTerm}
    \beta_{k + 1} - \tfrac{1}{2} H_{k + 1} r_{k + 1}^2
    \leq
    \PositivePart{\beta_{k + 1} - \tfrac{1}{2} H_{k + 1} r_{k + 1}^2}
    =
    (H_{k + 1} - H_k) D^2.
  \end{equation}
  Substituting
  \cref{%
    \LocalName{eq:UpperBoundOnDistanceTerm},%
    \LocalName{eq:UpperBoundOnBregmanTerm}%
  }
  into \cref{\LocalName{eq:PreliminaryRecursion}}, we get
  \begin{equation}
    \LocalLabel{eq:Recursion}
    F(x_{k + 1})
    +
    \tfrac{1}{2} H_{k + 1} \Norm{x - x_{k + 1}}^2
    \leq
    \zeta_k(x) + \tfrac{1}{2} H_k \Norm{x - x_k}^2
    +
    2 (H_{k + 1} - H_k) D^2.
  \end{equation}

  Let $k \geq 1$ be arbitrary.
  Summing up \cref{\LocalName{eq:Recursion}} for all indices $0 \leq k' \leq k -
  1$ and using the fact that $H_0 = 0$, we obtain
  \[
    \sum_{i = 1}^k F(x_i)
    \leq
    \sum_{i = 0}^{k - 1} \zeta_i(x) + 2 H_k D^2
    =
    k \Phi_k(x) + 2 H_k D^2,
  \]
  where the final identity follows from
  \cref{%
    \LocalName{eq:PartOfEstimatingFunction},%
    eq:EstimatingFunctionInGradientMethod%
  }
  (and we have dropped the nonnegative term $\frac{1}{2} H_k \Norm{x - x_k}^2$
  from the left-hand side).
  Since $x \in \EffectiveDomain \psi$ was arbitrary, this proves the second
  inequality in \cref{\MasterName{eq:Result}}.

  \ProofPart

  It remains to estimate the rate of growth of the coefficients $H_k$.

  Let $\nu \in \ClosedClosedInterval{0}{1}$ be arbitrary such that $H_{\nu} <
  +\infty$.
  From \cref{eq:UpperBoundOnBregmanDistanceViaHolderConstant} and the
  definitions of~$\beta_{k + 1}$ and~$r_{k + 1}$,
  we obtain $\beta_{k + 1} \leq \frac{L_{\nu}}{1 + \nu} r_{k + 1}^{1 + \nu}$ for
  any $k \geq 0$.
  Hence, according to \cref{\LocalName{eq:EquationForRegularizationParameter}},
  for all $k \geq 0$, we have the following bound:
  \[
    (H_{k + 1} - H_k) D^2
    \leq
    \PositivePart[\Big]{
      \frac{L_{\nu}}{1 + \nu} r_{k + 1}^{1 + \nu}
      -
      \frac{1}{2} H_{k + 1} r_{k + 1}^2
    }.
  \]
  Applying \cref{th:EstimatingGrowthRateOfRegularizationParameter} (with $\Omega
  \DefinedEqual D^2$, $M \DefinedEqual L_{\nu}$, $\gamma_k \DefinedEqual 1$), we
  get, for all $k \geq 1$,
  \[
    H_k
    \leq
    \Bigl[ \frac{1}{(1 + \nu) D^2} k \Bigr]^{(1 - \nu) / 2} L_{\nu}
    \leq
    \frac{L_{\nu}}{D^{1 - \nu}} k^{(1 - \nu) / 2}.
  \]
  Substituting this estimate into \cref{\MasterName{eq:Result}} and using the
  fact that $\nu \in \ClosedClosedInterval{0}{1}$ was arbitrary, we obtain the
  final inequality in \cref{\MasterName{eq:Result}}.
\end{proof}

%% file: Appendix/ProofForStochasticGradientMethod.tex
\section{Proof of Theorem~\ref{th:ConvergenceRateOfStochasticGradientMethod}}
\label{sec:ProofForStochasticGradientMethod}

\thConvergenceRateOfStochasticGradientMethod*

\begin{proof}
  \UsingNestedNamespace{ConvergenceRateOfStochasticGradientMethod}{Proof}

  \ProofPart

  We will show that
  \begin{equation}
    \LocalLabel{eq:Result}
    \Expectation[F(\bar{x}_k)] - F^*
    \leq
    \frac{2 \Expectation[H_k] D^2}{k}
    \leq
    \inf_{\nu \in \ClosedClosedInterval{0}{1}}
    \frac{8 L_{\nu} D^{1 + \nu}}{k^{(1 + \nu) / 2}}
    +
    \frac{4 \sigma D}{\sqrt{k}}.
  \end{equation}
  
  \ProofPart

  Let $k \geq 0$ be arbitrary.
  Let $\Map{\zeta_k}{\RealField^n}{\RealFieldPlusInfty}$ be the function
  \begin{equation}
    \LocalLabel{eq:StochasticLowerBoundFunction}
    \zeta_k(x)
    \DefinedEqual
    f(x_k) + \DualPairing{g_k}{x - x_k} + \psi(x).
  \end{equation}
  Note that, by definition, $g_k = g(x_k, \xi_k)$, where $\xi_k$ is independent
  of $x_k$.
  Therefore, in view of \cref{eq:StochasticGradientIsUnbiased,eq:MainProblem},
  in expectation, $\zeta_k$ is a global lower bound on the objective
  function~$F$: for all $x \in \EffectiveDomain \psi$, we have
  \begin{equation}
    \LocalLabel{eq:ExpectationOfStochasticLowerBoundFunction}
    \Expectation_{\xi_k}[\zeta_k(x)]
    =
    \Expectation_{\xi_k}[
      f(x_k)
      +
      \DualPairing{f'(x_k)}{x - x_k}
    ]
    +
    \psi(x)
    \leq
    f(x) + \psi(x)
    =
    F(x).
  \end{equation}

  Let $x \in \EffectiveDomain \psi$ be arbitrary.
  By the definition of $x_{k + 1}$
  (\cref{StochaticGradientMethod::step:ComputeNewPoint}) and
  \cref{th:OptimalityConditionForProximalPointStep},
  \begin{equation}
    \LocalLabel{eq:OptimalityCondition}
    \zeta_k(x) + \tfrac{1}{2} H_k \Norm{x - x_k}^2
    \geq
    \zeta_k(x_{k + 1})
    +
    \tfrac{1}{2} H_k \Norm{x_{k + 1} - x_k}^2
    +
    \tfrac{1}{2} H_k \Norm{x - x_{k + 1}}^2.
  \end{equation}
  According to
  \cref{\LocalName{eq:StochasticLowerBoundFunction},eq:MainProblem}, we have
  \begin{equation}
    \LocalLabel{eq:StochasticLowerBoundFunctionAtNewPoint}
    \zeta_k(x_{k + 1})
    =
    f(x_k) + \DualPairing{g_k}{x_{k + 1} - x_k} + \psi(x_{k + 1})
    =
    F(x_{k + 1}) - \beta_{k + 1},
  \end{equation}
  where $\beta_{k + 1} \DefinedEqual f(x_{k + 1}) - f(x_k) -
  \DualPairing{g_k}{x_{k + 1} - x_k}$.
  Using \cref{eq:StochasticGradientIsUnbiased}, we can estimate
  \begin{equation}
    \LocalLabel{eq:UpperBoundOnStochasticBregmanDistance}
    \beta_{k + 1}
    \leq
    \DualPairing{f'(x_{k + 1}) - g_k}{x_{k + 1} - x_k}
    =
    \DualPairing{g_{k + 1} - g_k}{x_{k + 1} - x_k}
    +
    \Delta_{k + 1},
  \end{equation}
  where $\Delta_{k + 1} \DefinedEqual \DualPairing{f'(x_{k + 1}) - g_{k +
  1}}{x_{k + 1} - x_k}$.
  Recall that $g_{k + 1} = g(x_{k + 1}, \xi_{k + 1})$ with $\xi_{k + 1}$ being
  independent of $x_k$ and $x_{k + 1}$.
  Therefore, according to \cref{eq:StochasticGradientIsUnbiased},
  \begin{equation}
    \LocalLabel{eq:ExpectedErrorInEstimatingStochasticBregmanDistance}
    \Expectation_{\xi_{k + 1}}[\Delta_{k + 1}] = 0.
  \end{equation}
  Substituting
  \cref{%
    \LocalName{eq:StochasticLowerBoundFunctionAtNewPoint},%
    \LocalName{eq:UpperBoundOnStochasticBregmanDistance}%
  }
  into \cref{\LocalName{eq:OptimalityCondition}}, we obtain
  \begin{align*}
    \zeta_k(x)
    +
    \tfrac{1}{2} H_k \Norm{x - x_k}^2
    &\geq
    F(x_{k + 1})
    -
    \DualPairing{g_{k + 1} - g_k}{x_{k + 1} - x_k}
    -
    \Delta_{k + 1}
    \\
    &\qquad
    +
    \tfrac{1}{2} H_k \Norm{x_{k + 1} - x_k}^2
    +
    \tfrac{1}{2} H_k \Norm{x - x_{k + 1}}^2
    \\
    &=
    F(x_{k + 1}) - \hat{\beta}_{k + 1} + \tfrac{1}{2} H_k r_{k + 1}^2
    +
    \tfrac{1}{2} H_k \Norm{x - x_{k + 1}}^2
    -
    \Delta_{k + 1},
  \end{align*}
  where the last identity is due to the definitions of~$r_{k + 1}$
  and~$\beta_{k + 1}$.
  Rearranging, we get
  \begin{equation}
    \LocalLabel{eq:PreliminaryRecursion}
    \begin{aligned}
      F(x_{k + 1})
      +
      \tfrac{1}{2} H_{k + 1} \Norm{x - x_{k + 1}}^2
      &\leq
      \zeta_k(x)
      +
      \tfrac{1}{2} H_k \Norm{x - x_k}^2
      +
      [\hat{\beta}_{k + 1} - \tfrac{1}{2} H_{k + 1} r_{k + 1}^2]
      \\
      &\qquad
      +
      \tfrac{1}{2} (H_{k + 1} - H_k)
      (\Norm{x - x_{k + 1}}^2 + r_{k + 1}^2)
      +
      \Delta_{k + 1}.
    \end{aligned}
  \end{equation}
  By construction, $H_k \leq H_{k + 1}$
  (\cref{StochaticGradientMethod::step:UpdateRegularizationParameter}).
  Also, in view of \cref{as:BoundedFeasibleSet} (and the fact that $x_i \in
  \EffectiveDomain \psi$ for all $i \geq 0$), $\Norm{x - x_{k + 1}} \leq D$ and
  $r_{k + 1} \leq D$.
  Therefore,
  \begin{equation}
    \LocalLabel{eq:UpperBoundOnDistanceTerm}
    \tfrac{1}{2} (H_{k + 1} - H_k)
    (\Norm{x - x_{k + 1}}^2 + r_{k + 1}^2)
    \leq
    (H_{k + 1} - H_k) D^2.
  \end{equation}
  Further, by the definition of $H_{k + 1}$ at
  \cref{StochaticGradientMethod::step:UpdateRegularizationParameter}, it
  satisfies the following equation (see
  \cref{th:SolvingEquationForRegularizationParameter}):
  \begin{equation}
    \LocalLabel{eq:EquationForRegularizationParameter}
    (H_{k + 1} - H_k) D^2
    =
    \PositivePart{
      \hat{\beta}_{k + 1} - \tfrac{1}{2} H_{k + 1} r_{k + 1}^2
    }.
  \end{equation}
  Substituting
  \cref{%
    \LocalName{eq:UpperBoundOnDistanceTerm},%
    \LocalName{eq:EquationForRegularizationParameter}%
  }
  into \cref{\LocalName{eq:PreliminaryRecursion}}, we obtain
  \begin{equation}
    \LocalLabel{eq:Recursion}
    F(x_{k + 1})
    +
    \tfrac{1}{2} H_{k + 1} \Norm{x - x_{k + 1}}^2
    \leq
    \zeta_k(x)
    +
    \tfrac{1}{2} H_k \Norm{x - x_k}^2
    +
    2 (H_{k + 1} - H_k) D^2
    +
    \Delta_{k + 1}.
  \end{equation}

  Let $k \geq 1$ be arbitrary.
  Summing up \cref{\LocalName{eq:Recursion}} for all indices
  $0 \leq k' \leq k - 1$ and using the fact that $H_0 = 0$, we get
  \[
    \sum_{i = 1}^k F(x_i)
    \leq
    \sum_{i = 0}^{k - 1} \zeta_i(x)
    +
    2 H_k D^2
    +
    \sum_{i = 1}^k \Delta_i.
  \]
  Hence, by the convexity of $F$ and the definition of~$\bar{x}_k$,
  \[
    F(\bar{x}_k)
    \leq
    \frac{1}{k} \sum_{i = 1}^k F(x_i)
    \leq
    \frac{1}{k} \sum_{i = 0}^{k - 1} \zeta_i(x)
    +
    \frac{2 H_k D^2}{k}
    +
    \frac{1}{k} \sum_{i = 1}^k \Delta_i.
  \]
  Passing to expectations and taking into account
  \cref{%
    \LocalName{eq:ExpectationOfStochasticLowerBoundFunction},%
    \LocalName{eq:ExpectedErrorInEstimatingStochasticBregmanDistance}%
  },
  we conclude that
  \[
    \Expectation[F(\bar{x}_k)]
    \leq
    F(x)
    +
    \frac{2 \Expectation[H_k] D^2}{k}.
  \]
  This proves the first inequality in \cref{\LocalName{eq:Result}} since
  $x \in \EffectiveDomain \psi$ was arbitrary.

  \ProofPart

  Let us estimate the rate of growth of the coefficients $H_k$.

  For each $k \geq 0$, denote
  \begin{equation}
    \LocalLabel{eq:GradientEstimateError}
    \delta_k
    \DefinedEqual
    g_k - f'(x_k)
    =
    g(x_k, \xi_k) - f'(x_k).
  \end{equation}
  Note that, according to
  \cref{%
    eq:StochasticGradientIsUnbiased,%
    eq:VarianceOfStochasticGradient%
  },
  we have
  \begin{equation}
    \LocalLabel{eq:ExpectationAndVarianceOfGradientEstimateError}
    \Expectation_{\xi_k}[\delta_k] = 0,
    \qquad
    \Expectation_{\xi_k}[\DualNorm{\delta_k}^2] \leq \sigma^2.
  \end{equation}

  Let $\nu \in \ClosedClosedInterval{0}{1}$ be arbitrary such that
  $L_{\nu} < +\infty$.
  Let $k \geq 0$ be arbitrary.
  By the definitions of $\hat{\beta}_{k + 1}$ and $r_{k + 1}$, and by
  \cref{\LocalName{eq:GradientEstimateError}}, we have
  \begin{equation}
    \LocalLabel{eq:UpperBoundOnStochasticSymmetrizedBregmanDistance}
    \begin{aligned}
      \hat{\beta}_{k + 1}
      &=
      \DualPairing{g_{k + 1} - g_k}{x_{k + 1} - x_k}
      \\
      &=
      \DualPairing{f'(x_{k + 1}) - f'(x_k)}{x_{k + 1} - x_k}
      +
      \DualPairing{\delta_{k + 1} - \delta_k}{x_{k + 1} - x_k}
      \\
      &\leq
      \DualNorm{f'(x_{k + 1}) - f'(x_k)} r_{k + 1}
      +
      \DualNorm{\delta_{k + 1} - \delta_k} r_{k + 1}
      \leq
      L_{\nu} r_{k + 1}^{1 + \nu}
      +
      \sigma_{k + 1} r_{k + 1},
    \end{aligned}
  \end{equation}
  where $\sigma_{k + 1} \DefinedEqual \DualNorm{\delta_{k + 1} - \delta_k}$, the
  first inequality is the Cauchy--Schwartz, and the final one is due to
  \cref{eq:HolderConstant,eq:StochasticGradientIsUnbiased}.
  Recall that $\xi_{k + 1}$ is independent of $x_k$ and $\xi_k$.
  Hence, it is also independent of $\delta_k$ (see
  \cref{\LocalName{eq:GradientEstimateError}}).
  Therefore, according to
  \cref{%
    eq:DualNorm,%
    \LocalName{eq:ExpectationAndVarianceOfGradientEstimateError}%
  },
  we have
  \begin{equation}
    \LocalLabel{eq:SecondMomentOfGradientDifferenceError}
    \begin{aligned}
      \Expectation_{\xi_k, \xi_{k + 1}}[\sigma_{k + 1}^2]
      &=
      \Expectation_{\xi_k, \xi_{k + 1}}[
        \DualNorm{\delta_{k + 1}}^2
        +
        \DualNorm{\delta_k}^2
        +
        \DualPairing{\delta_{k + 1}}{B^{-1} \delta_k}
      ]
      \\
      &=
      \Expectation_{\xi_k} \bigl[
        \Expectation_{\xi_{k + 1}}[\DualNorm{\delta_{k + 1}}^2]
        +
        \DualNorm{\delta_k}^2
      \bigr]
      \leq
      \Expectation_{\xi_k}[
        \sigma^2 + \DualNorm{\delta_k}^2
      ]
      \leq
      2 \sigma^2.
    \end{aligned}
  \end{equation}
  In particular, the same inequality holds for the full expectation.

  Substituting
  \cref{\LocalName{eq:UpperBoundOnStochasticSymmetrizedBregmanDistance}} into
  \cref{\LocalName{eq:EquationForRegularizationParameter}} (using the
  monotonicity of $\PositivePart{\cdot}$), we get, for any $k \geq 0$,
  \[
    (H_{k + 1} - H_k) D^2
    \leq
    \PositivePart{
      L_{\nu} r_{k + 1}^{1 + \nu}
      +
      \sigma_{k + 1} r_{k + 1}
      -
      \tfrac{1}{2} H_{k + 1} r_{k + 1}^2
    }.
  \]
  Applying
  \cref{th:EstimatingGrowthRateOfRegularizationParameterInStochasticCase}
  (with
  $\Omega \DefinedEqual D^2$, $L \DefinedEqual L_{\nu}$,
  $\alpha_k \DefinedEqual 1$, $\gamma_k \DefinedEqual \sigma_k$),
  we obtain, for all $k \geq 1$,
  \[
    H_k
    \leq
    [2 (1 + \nu)]^{(1 + \nu) / 2}
    \frac{L_{\nu}}{D^{1 - \nu}}
    k^{(1 - \nu) / 2}
    +
    \Bigl( \frac{2}{D^2} \sum_{i = 1}^k \sigma_i^2 \Bigr)^{1 / 2}.
  \]
  By Jensen's inequality
  $\Expectation[X^{1 / 2}] \leq (\Expectation[X])^{1 / 2}$
  and \cref{\LocalName{eq:SecondMomentOfGradientDifferenceError}},
  for all $k \geq 1$, it holds
  \[
    \Expectation\biggl[
      \Bigl( \frac{2}{D^2} \sum_{i = 1}^k \sigma_i^2 \Bigr)^{1 / 2}
    \biggr]
    \leq
    \Bigl(
      \frac{2}{D^2} \sum_{i = 1}^k \Expectation[\sigma_i^2]
    \Bigr)^{1 / 2}
    \leq
    \sqrt{\frac{2}{D^2} (2 \sigma^2) k}
    =
    \frac{2 \sigma}{D} \sqrt{k}.
  \]
  Thus, for all $k \geq 1$,
  \[
    \Expectation[H_k]
    \leq
    [2 (1 + \nu)]^{(1 + \nu) / 2}
    \frac{L_{\nu}}{D^{1 - \nu}}
    k^{(1 - \nu) / 2}
    +
    \frac{2 \sigma}{D} \sqrt{k}
    \leq
    \frac{4 L_{\nu}}{D^{1 - \nu}} k^{(1 - \nu) / 2}
    +
    \frac{2 \sigma}{D} \sqrt{k}.
  \]
  Substituting this estimate into \cref{\LocalName{eq:Result}} and taking into
  account the fact that $\nu \in \ClosedClosedInterval{0}{1}$ was arbitrary, we
  obtain the final inequality in \cref{\LocalName{eq:Result}}.
\end{proof}

%% file: Appendix/ProofForStochasticFastGradientMethod.tex
\section{Proof of Theorem~\ref{th:ConvergenceRateOfStochasticFastGradientMethod}}
\label{sec:ProofForStochasticFastGradientMethod}

\thConvergenceRateOfStochasticFastGradientMethod*

\begin{proof}
  \UsingNestedNamespace{ConvergenceRateOfStochasticFastGradientMethod}{Proof}

  \ProofPart

  We will show that
  \begin{equation}
    \LocalLabel{eq:Result}
    \Expectation[F(x_k)] - F^*
    \leq
    \frac{2 \Expectation[H_k] D^2}{A_k}
    \leq
    \frac{4 \Expectation[H_k] D^2}{k (k + 1)}
    \leq
    \inf_{\nu \in \ClosedClosedInterval{0}{1}}
    \frac{32 L_{\nu} D^{1 + \nu}}{k^{(1 + 3 \nu) / 2}}
    +
    \frac{8 \sigma D}{\sqrt{3 k}}.
  \end{equation}

  \ProofPart

  For each $k \geq 0$, denote
  \begin{gather}
    \LocalLabel{eq:GradientApproximationErrorAtIntermediatePoint}
    \delta_k^y
    \DefinedEqual
    g_k^y - f'(y_k)
    \equiv
    g(y_k, \xi_k^y) - f'(y_k),
    \\
    \LocalLabel{eq:GradientApproximationErrorAtNewPoint}
    \delta_{k + 1}^x
    \DefinedEqual
    g_{k + 1}^x - f'(x_{k + 1})
    \equiv
    g(x_{k + 1}, \xi_{k + 1}^x) - f'(x_{k + 1}),
  \end{gather}
  where $\xi_k^y$, $\xi_{k + 1}^x$, $k = 0, 1, \ldots$ are independent copies of
  the oracle's randomness.
  
  Note that $\xi_k^y$ (resp., $\xi_{k + 1}^y$) is generated after (and
  independently) of~$y_k$ (resp., $x_{k + 1}$).
  Therefore, according to
  \cref{eq:StochasticGradientIsUnbiased,eq:VarianceOfStochasticGradient},
  \begin{gather}
    \LocalLabel{eq:ExpectationAndVarianceOfGradientApproximationErrorAtIntermediatePoint}
    \Expectation_{\xi_k^y} [\delta_k^y] = 0,
    \qquad
    \Expectation_{\xi_k^y} [\DualNorm{\delta_k^y}^2] \leq \sigma^2,
    \\
    \LocalLabel{eq:ExpectationAndVarianceOfGradientApproximationErrorAtNewPoint}
    \Expectation_{\xi_{k + 1}^x} [\delta_{k + 1}^x] = 0,
    \qquad
    \Expectation_{\xi_{k + 1}^x} [\DualNorm{\delta_{k + 1}^x}^2] \leq \sigma^2.
  \end{gather}

  \ProofPart

  Let $k \geq 0$ be arbitrary.
  Let $\Map{\zeta_k}{\RealField^n}{\RealFieldPlusInfty}$ be the global lower
  bound on the objective function~$F$ obtained by linearizing~$f$ at~$y_k$:
  \begin{equation}
    \LocalLabel{eq:LowerBoundFunction}
    \zeta_k(x)
    \DefinedEqual
    f(y_k) + \DualPairing{f'(y_k)}{x - y_k} + \psi(x)
    \quad
    (\leq F(x)),
  \end{equation}
  and let $\Map{\hat{\zeta}_k}{\RealField^n}{\RealFieldPlusInfty}$ be its
  stochastic approximation:
  \begin{equation}
    \LocalLabel{eq:ApproximationOfLowerBoundFunction}
    \hat{\zeta}_k(x)
    \DefinedEqual
    f(y_k) + \DualPairing{g_k^y}{x - y_k} + \psi(x)
    =
    \zeta_k(x) + \Delta_k^y(x),
  \end{equation}
  where
  \begin{equation}
    \LocalLabel{eq:ApproximationErrorForLowerBoundFunction}
    \Delta_k^y(x) \DefinedEqual \DualPairing{\delta_k^y}{x - y_k}.
  \end{equation}
  Recall that $\xi_k^y$ is independent of~$y_k$.
  Therefore, for any (possibly random variable) $x \in \RealField^n$
  that is also independent of~$\xi_k^y$, we have
  \begin{equation}
    \LocalLabel{eq:ExpectedApproximationErrorForLowerBoundFunction}
    \Expectation_{\xi_k^y} [\Delta_k^y(x)] = 0.
  \end{equation}

  Let $x \in \EffectiveDomain \psi$ be an arbitrary (deterministic) point.
  Applying \cref{th:OptimalityConditionForProximalPointStep}
  to the definition of~$v_{k + 1}$ at
  \cref{StochasticFastGradientMethod::step:UpdateProxCenter}
  and taking into account
  \cref{\LocalName{eq:ApproximationOfLowerBoundFunction}},
  we obtain
  \begin{equation}
    \LocalLabel{eq:OptimalityCondition}
    \begin{aligned}
      a_{k + 1} \hat{\zeta}_k(x) + \tfrac{1}{2} H_k \Norm{x - v_k}^2
      &\geq
      a_{k + 1} \hat{\zeta}_k(v_{k + 1})
      +
      \tfrac{1}{2} H_k \Norm{v_{k + 1} - v_k}^2
      +
      \tfrac{1}{2} H_k \Norm{x - v_{k + 1}}^2
      \\
      &=
      a_{k + 1} \hat{\zeta}_k(v_{k + 1})
      +
      \tfrac{1}{2} H_k \Norm{x - v_{k + 1}}^2
      +
      \tfrac{1}{2} H_k r_{k + 1}^2,
    \end{aligned}
  \end{equation}
  where the last identity follows from the definition of~$r_{k + 1}$.

  In view of
  \cref{%
    \LocalName{eq:LowerBoundFunction},%
    \LocalName{eq:ApproximationOfLowerBoundFunction}%
  },
  we have
  \begin{equation}
    \LocalLabel{eq:PreliminaryLowerBoundOnCombinationOfFunctonValues}
    \begin{aligned}
      \hspace{2em}&\hspace{-2em}
      A_k F(x_k) + a_{k + 1} \hat{\zeta}_k(v_{k + 1})
      \geq
      A_k \zeta_k(x_k) + a_{k + 1} \hat{\zeta}_k(v_{k + 1})
      \\
      &=
      A_k \hat{\zeta}_k(x_k)
      +
      a_{k + 1} \hat{\zeta}_k(v_{k + 1})
      -
      A_k \Delta_k^y(x_k)
      \geq
      A_{k + 1} \hat{\zeta}_k(x_{k + 1}) - A_k \Delta_k^y(x_k),
    \end{aligned}
  \end{equation}
  where the last inequality follows from the convexity of~$\hat{\zeta}_k$
  and the definitions of~$x_{k + 1}$ and~$A_{k + 1}$ at
  \cref{%
    StochasticFastGradientMethod::step:ComputeNewPoint,%
    StochasticFastGradientMethod::step:ChooseScalingCoefficient%
  },
  respectively.
  According to
  \cref{%
    \LocalName{eq:ApproximationOfLowerBoundFunction},%
    eq:MainProblem%
  },
  \begin{equation}
    \LocalLabel{eq:ApproximationOfLowerBoundFunctionAtNewPoint}
    \hat{\zeta}_k(x_{k + 1})
    =
    f(y_k) + \DualPairing{g_k^y}{x_{k + 1} - y_k} + \psi(x_{k + 1})
    =
    F(x_{k + 1}) - \beta_{k + 1},
  \end{equation}
  where
  $
    \beta_{k + 1}
    \DefinedEqual
    f(x_{k + 1}) - f(y_k) - \DualPairing{g_k^y}{x_{k + 1} - y_k}
  $.
  Using \cref{eq:StochasticGradientIsUnbiased} and the definition
  of~$\hat{\beta}_{k + 1}$, we can estimate
  \begin{equation}
    \LocalLabel{eq:UpperBoundOnBregmanDistanceApproximation}
    \beta_{k + 1}
    \leq
    \DualPairing{f'(x_{k + 1}) - g_k^y}{x_{k + 1} - y_k}
    =
    \hat{\beta}_{k + 1} + \Delta_{k + 1}^x,
  \end{equation}
  where
  $
    \Delta_{k + 1}^x
    \DefinedEqual
    \DualPairing{\delta_{k + 1}^x}{y_k - x_{k + 1}}
  $
  (see \cref{\LocalName{eq:GradientApproximationErrorAtNewPoint}}).
  Note that $\xi_{k + 1}^x$ is generated after (and independently)
  of~$x_{k + 1}$ and~$y_k$.
  Hence, in view of
  \cref{\LocalName{eq:ExpectationAndVarianceOfGradientApproximationErrorAtNewPoint}},
  \begin{equation}
    \LocalLabel{eq:ExpectedErrorOfStochasticBregmanDistanceApproximation}
    \Expectation_{\xi_{k + 1}^x} [\Delta_{k + 1}^x] = 0.
  \end{equation}
  Putting together
  \cref{%
    \LocalName{eq:PreliminaryLowerBoundOnCombinationOfFunctonValues},%
    \LocalName{eq:ApproximationOfLowerBoundFunctionAtNewPoint},%
    \LocalName{eq:UpperBoundOnBregmanDistanceApproximation}%
  }
  we obtain
  \begin{align*}
    A_k F(x_k) + a_{k + 1} \hat{\zeta}_k(v_{k + 1})
    &\geq
    A_{k + 1} [F(x_{k + 1}) - \beta_{k + 1}] - A_k \Delta_k^y(x_k)
    \\
    &\geq
    A_{k + 1} [F(x_{k + 1}) - \hat{\beta}_{k + 1}]
    -
    A_k \Delta_k^y(x_k)
    -
    A_{k + 1} \Delta_{k + 1}^x.
  \end{align*}

  Combining the above inequality with
  \cref{\LocalName{eq:OptimalityCondition}}, we get
  \begin{align*}
    \hspace{2em}&\hspace{-2em}
    A_k F(x_k) + a_{k + 1} \hat{\zeta}_k(x)
    +
    \tfrac{1}{2} H_k \Norm{x - v_k}^2
    \\
    &\geq
    A_{k + 1} [F(x_{k + 1}) - \hat{\beta}_{k + 1}]
    +
    \tfrac{1}{2} H_k r_{k + 1}^2
    +
    \tfrac{1}{2} H_k \Norm{x - v_{k + 1}}^2
    -
    A_k \Delta_k^y(x_k)
    -
    A_{k + 1} \Delta_{k + 1}^x.
  \end{align*}
  After rearranging, we can write
  \begin{align*}
    \hspace{2em}&\hspace{-2em}
    A_{k + 1} F(x_{k + 1}) + \tfrac{1}{2} H_{k + 1} \Norm{x - v_{k + 1}}^2
    \\
    &\leq
    A_k F(x_k) + \tfrac{1}{2} H_k \Norm{x - v_k}^2
    +
    a_{k + 1} \hat{\zeta}_k(x)
    +
    [
      A_{k + 1} \hat{\beta}_{k + 1}
      -
      \tfrac{1}{2} H_{k + 1} r_{k + 1}^2
    ]
    \\
    &\qquad
    +
    \tfrac{1}{2} (H_{k + 1} - H_k) [\Norm{x - v_{k + 1}}^2 + r_{k + 1}^2]
    +
    A_k \Delta_k^y(x_k)
    +
    A_{k + 1} \Delta_{k + 1}^x.
  \end{align*}
  Note that, by construction, $H_k \leq H_{k + 1}$
  (\cref{StochasticFastGradientMethod::step:UpdateRegularizationParameter}).
  Further, in view of \cref{as:BoundedFeasibleSet}
  (and the fact that $v_i \in \EffectiveDomain \psi$ for all $i \geq 0$),
  we have $\Norm{x - v_{k + 1}} \leq D$ and $r_{k + 1} \leq D$.
  Hence,
  \[
    \tfrac{1}{2} (H_{k + 1} - H_k) [\Norm{x - v_{k + 1}}^2 + r_{k + 1}^2]
    \leq
    (H_{k + 1} - H_k) D^2.
  \]
  On the other hand, from the definition of~$H_{k + 1}$ at
  \cref{StochasticFastGradientMethod::step:UpdateRegularizationParameter}
  and \cref{th:SolvingEquationForRegularizationParameter}
  (with
  $\beta \DefinedEqual A_{k + 1} \hat{\beta}_{k + 1}$,
  $\rho \DefinedEqual \frac{1}{2} r_{k + 1}^2$,
  $\Omega \DefinedEqual D^2$),
  it follows that
  \begin{equation}
    \LocalLabel{eq:EquationForRegularizationParameter}
    (H_{k + 1} - H_k) D^2
    =
    \PositivePart{
      A_{k + 1} \hat{\beta}_{k + 1} - \tfrac{1}{2} H_{k + 1} r_{k + 1}^2
    }.
  \end{equation}
  Combining the above three displays, we obtain
  \begin{equation}
    \LocalLabel{eq:Recursion}
    \begin{multlined}
      A_{k + 1} F(x_{k + 1}) + \tfrac{1}{2} H_{k + 1} \Norm{x - v_{k + 1}}^2
      \\
      \leq
      A_k F(x_k) + \tfrac{1}{2} H_k \Norm{x - v_k}^2
      +
      a_{k + 1} \hat{\zeta}_k(x)
      +
      2 (H_{k + 1} - H_k) D^2
      \\
      +
      A_k \Delta_k^y(x_k)
      +
      A_{k + 1} \Delta_{k + 1}^x.
    \end{multlined}
  \end{equation}
  Note that this inequality is valid for any $k \geq 0$.

  Let $k \geq 1$ be arbitrary.
  Summing up \cref{\LocalName{eq:Recursion}} for all
  indices~$0 \leq k' \leq k - 1$ and taking into account that
  $H_0 = A_0 = 0$ (by definition), we get
  \[
    A_k F(x_k)
    \leq
    \sum_{i = 0}^{k - 1} a_{i + 1} \hat{\zeta}_i(x)
    +
    2 H_k D^2
    +
    \sum_{i = 0}^{k - 1} (A_i \Delta_i^y(x_i) + A_{i + 1} \Delta_{i + 1}^x),
  \]
  where we have additionally dropped the nonnegative term
  $\frac{1}{2} H_k \Norm{x - v_k}^2$ from the left-hand side.
  Combining this with
  \cref{%
    \LocalName{eq:ApproximationOfLowerBoundFunction},%
    \LocalName{eq:LowerBoundFunction}%
  },
  we obtain
  \begin{align*}
    A_k F(x_k)
    &\leq
    \sum_{i = 0}^{k - 1} a_{i + 1} [\zeta_i(x) + \Delta_i^y(x)]
    +
    2 H_k D^2
    +
    \sum_{i = 0}^{k - 1} (A_i \Delta_i^y(x_i) + A_{i + 1} \Delta_{i + 1}^x)
    \\
    &\leq
    A_k F(x) + 2 H_k D^2
    +
    \sum_{i = 0}^{k - 1}
    (
      a_{i + 1} \Delta_i^y(x)
      +
      A_i \Delta_i^y(x_i)
      +
      A_{i + 1} \Delta_{i + 1}^x
    ),
  \end{align*}
  where we have used the fact that $A_k = \sum_{i = 1}^k a_i$
  (see \cref{StochasticFastGradientMethod::step:ChooseScalingCoefficient}).

  Observe that, by definitions at
  \cref{StochasticFastGradientMethod::step:ChooseScalingCoefficient},
  the coefficients~$a_i$ and~$A_i$ are deterministic for each $i \geq 0$.
  Also, recall that $x$ is assumed to be deterministic as well.
  Therefore, passing to expectations in the above inequality, we get
  \[
    A_k \Expectation [F(x_k)]
    \leq
    A_k F(x)
    +
    2 \Expectation [H_k] D^2
    +
    \sum_{i = 0}^{k - 1}
    (
      a_{i + 1} \Expectation [\Delta_i^y(x)]
      +
      A_i \Expectation [\Delta_i^y(x_i)]
      +
      A_{i + 1} \Expectation [\Delta_{i + 1}^x]
    ).
  \]
  Note that, for any $i \geq 0$, the random variable~$\xi_i^y$ is generated
  after $x_i$, and hence they are independent.
  Therefore, according to
  \cref{%
    \LocalName{eq:ApproximationErrorForLowerBoundFunction},%
    \LocalName{eq:ExpectedErrorOfStochasticBregmanDistanceApproximation}%
  },
  for each $i \geq 0$, we have
  $
    \Expectation [\Delta_i^y(x)]
    =
    \Expectation [\Delta_i^y(x_i)]
    =
    \Expectation [\Delta_{i + 1}^x]
    =
    0
  $.
  Thus, the above display reads
  \[
    A_k \Expectation [F(x_k)] \leq A_k F(x) + 2 \Expectation [H_k] D^2.
  \]
  This proves the first inequality in \cref{\LocalName{eq:Result}}
  since $x \in \EffectiveDomain \psi$ was arbitrary.

  \ProofPart

  From the definitions at
  \cref{StochasticFastGradientMethod::step:ChooseScalingCoefficient}
  and the fact that $A_0 = 0$, it follows that
  \begin{equation}
    \LocalLabel{eq:ExplicitFormulaForScalingCoefficient}
    A_k
    =
    \sum_{i = 1}^k a_i
    =
    \sum_{i = 1}^k i
    =
    \tfrac{1}{2} k (k + 1)
    \quad
    (\geq \tfrac{1}{2} k^2)
  \end{equation}
  for any $k \geq 1$.
  This proves the second inequality in \cref{\LocalName{eq:Result}}.

  \ProofPart

  To prove the final inequality in \cref{\LocalName{eq:Result}}, it remains to
  estimate the expected growth rate of regularization parameters~$H_k$.

  Let $\nu \in \ClosedClosedInterval{0}{1}$ be arbitrary such that
  $L_{\nu} < +\infty$.
  Let $k \geq 0$ be arbitrary.
  According to the definition of~$\hat{\beta}_{k + 1}$
  and
  \cref{%
    \LocalName{eq:GradientApproximationErrorAtIntermediatePoint},%
    \LocalName{eq:GradientApproximationErrorAtNewPoint}%
  },
  we have
  \begin{equation}
    \LocalLabel{eq:UpperBoundOnStochasticSymmetrizedBregmanDistance}
    \begin{aligned}
      \hat{\beta}_{k + 1}
      &=
      \DualPairing{g_{k + 1}^x - g_k^y}{x_{k + 1} - y_k}
      \\
      &=
      \DualPairing{f'(x_{k + 1}) - f'(y_k)}{x_{k + 1} - y_k}
      +
      \DualPairing{\delta_{k + 1}^x - \delta_k^y}{x_{k + 1} - y_k}
      \\
      &\leq
      \DualNorm{f'(x_{k + 1}) - f'(y_k)} \Norm{x_{k + 1} - y_k}
      +
      \sigma_{k + 1} \Norm{x_{k + 1} - y_k}
      \\
      &\leq
      L_{\nu} \Norm{x_{k + 1} - y_k}^{1 + \nu}
      +
      \sigma_{k + 1} \Norm{x_{k + 1} - y_k},
    \end{aligned}
  \end{equation}
  where
  $\sigma_{k + 1} \DefinedEqual \DualNorm{\delta_{k + 1}^x - \delta_k^y}$,
  the first inequality is the Cauchy--Schwarz
  inequality, and the final one is due to
  \cref{eq:HolderConstant,eq:StochasticGradientIsUnbiased}.
  Recall that $\delta_k^y$ is a function of~$y_k$ and~$\xi_k^y$
  (see \cref{\LocalName{eq:GradientApproximationErrorAtIntermediatePoint}}),
  and $\xi_{k + 1}^x$ is generated after $y_k$ and~$\xi_k^y$.
  Therefore, $\xi_{k + 1}^x$ is independent of~$\delta_k^y$.
  Hence, according to
  \cref{%
    eq:DualNorm,%
    \LocalName{eq:ExpectationAndVarianceOfGradientApproximationErrorAtIntermediatePoint},%
    \LocalName{eq:ExpectationAndVarianceOfGradientApproximationErrorAtNewPoint}%
  },
  \begin{equation}
    \LocalLabel{eq:SecondMomentOfGradientDifferenceApproximationError}
    \begin{aligned}
      \Expectation_{\xi_k^y, \xi_{k + 1}^x} [\sigma_{k + 1}^2]
      &=
      \Expectation_{\xi_k^y, \xi_{k + 1}^x} [
        \DualNorm{\delta_{k + 1}^x}^2
        +
        \DualNorm{\delta_k^y}^2
        +
        2 \DualPairing{\delta_{k + 1}^x}{B^{-1} \delta_k^y}
      ]
      \\
      &=
      \Expectation_{\xi_k^y} \bigl[
        \Expectation_{\xi_{k + 1}^x} [\DualNorm{\delta_{k + 1}^x}^2]
        +
        \DualNorm{\delta_k^y}^2
      \bigr]
      \leq
      \Expectation_{\xi_k^y} [\sigma^2 + \DualNorm{\delta_k^y}^2]
      \leq
      2 \sigma^2.
    \end{aligned}
  \end{equation}
  Further, by the definitions of~$y_k$ and~$x_{k + 1}$ at
  \cref{%
    StochasticFastGradientMethod::step:ComputeIntermediatePoint,%
    StochasticFastGradientMethod::step:ComputeNewPoint%
  },
  $x_{k + 1} - y_k = \frac{a_{k + 1}}{A_{k + 1}} (v_{k + 1} - v_k)$,
  which means that
  $\Norm{x_{k + 1} - y_k} = \frac{a_{k + 1}}{A_{k + 1}} r_{k + 1}$.
  Substituting this into
  \cref{\LocalName{eq:UpperBoundOnStochasticSymmetrizedBregmanDistance}}
  and using the definition of~$a_{k + 1}$ at
  \cref{StochasticFastGradientMethod::step:ChooseScalingCoefficient}
  together with \cref{\LocalName{eq:ExplicitFormulaForScalingCoefficient}},
  we obtain
  \begin{align*}
    A_{k + 1} \hat{\beta}_{k + 1}
    &\leq
    A_{k + 1}
    \Bigl[
      L_{\nu} \Bigl( \frac{a_{k + 1}}{A_{k + 1}} r_{k + 1} \Bigr)^{1 + \nu}
      +
      \sigma_{k + 1} \frac{a_{k + 1}}{A_{k + 1}} r_{k + 1}
    \Bigr]
    \\
    &=
    L_{\nu} \frac{a_{k + 1}^{1 + \nu}}{A_{k + 1}^\nu} r_{k + 1}^{1 + \nu}
    +
    a_{k + 1}\sigma_{k + 1} r_{k + 1}
    \\
    &\leq
    L_{\nu}
    \frac{(k + 1)^{1 + \nu}}{[(\frac{1}{2} (k + 1)^2]^{\nu}}
    r_{k + 1}^{1 + \nu}
    +
    (k + 1) \sigma_{k + 1} r_{k + 1}
    \\
    &=
    2^{\nu} L_{\nu} (k + 1)^{1 - \nu} r_{k + 1}^{1 + \nu}
    +
    (k + 1) \sigma_{k + 1} r_{k + 1}.
  \end{align*}
  Combining the above inequality with
  \cref{\LocalName{eq:EquationForRegularizationParameter}}
  (using the monotonicity of~$\PositivePart{\cdot}$), we come to the following
  recurrence relation:
  \[
    (H_{k + 1} - H_k) D^2
    \leq
    \PositivePart{
      2^{\nu} L_{\nu} (k + 1)^{1 - \nu} r_{k + 1}^{1 + \nu}
      +
      (k + 1) \sigma_{k + 1} r_{k + 1}
      -
      \tfrac{1}{2} H_{k + 1} r_{k + 1}^2
    },
  \]
  which is valid for any $k \geq 0$.

  Let $k \geq 1$ be arbitrary.
  Applying
  \cref{th:EstimatingGrowthRateOfRegularizationParameterInStochasticCase}
  (with
  $\Omega \DefinedEqual D^2$,
  $L \DefinedEqual 2^{\nu} L_{\nu}$,
  $\alpha_k \DefinedEqual k$
  and~$\gamma_k \DefinedEqual k \sigma_k$),
  we conclude that
  \begin{align*}
    H_k
    &\leq
    [2 (1 + \nu)]^{(1 + \nu) / 2}
    2^{\nu} L_{\nu}
    \Bigl( \frac{1}{D^2} \sum_{i = 1}^k i^2 \Bigr)^{(1 - \nu) / 2}
    +
    \Bigl( \frac{2}{D^2} \sum_{i = 1}^k (i \sigma_i)^2 \Bigr)^{1 / 2}
    \\
    &=
    2^{(1 + 3 \nu) / 2} (1 + \nu)^{(1 + \nu) / 2}
    \frac{L_{\nu}}{D^{1 - \nu}}
    \Bigl( \sum_{i = 1}^k i^2 \Bigr)^{(1 - \nu) / 2}
    +
    \Bigl( \frac{2}{D^2} \sum_{i = 1}^k i^2 \sigma_i^2 \Bigr)^{1 / 2}.
  \end{align*}
  Note that, by Jensen's inequality
  $\Expectation [X^{1 / 2}] \leq (\Expectation [X])^{1 / 2}$
  and \cref{\LocalName{eq:SecondMomentOfGradientDifferenceApproximationError}},
  \[
    \Expectation \biggl[
      \Bigl( \frac{2}{D^2} \sum_{i = 1}^k i^2 \sigma_i^2 \Bigr)^{1 / 2}
    \biggr]
    \leq
    \Bigl(
      \frac{2}{D^2} \sum_{i = 1}^k i^2 \Expectation [\sigma_i^2]
    \Bigr)^{1 / 2}
    \leq
    \frac{2 \sigma}{D} \Bigl( \sum_{i = 1}^k i^2 \Bigr)^{1 / 2}.
  \]
  Thus,
  \begin{align*}
    \Expectation [H_k]
    &\leq
    2^{(1 + 3 \nu) / 2} (1 + \nu)^{(1 + \nu) / 2}
    \frac{L_{\nu}}{D^{1 - \nu}}
    \Bigl( \sum_{i = 1}^k i^2 \Bigr)^{(1 - \nu) / 2}
    +
    \frac{2 \sigma}{D} \Bigl( \sum_{i = 1}^k i^2 \Bigr)^{1 / 2}
    \\
    &\leq
    2^{(1 + 3 \nu) / 2} (1 + \nu)^{(1 + \nu) / 2}
    \frac{L_{\nu}}{D^{1 - \nu}}
    \Bigl( \frac{1}{3} k (k + 1)^2 \Bigr)^{(1 - \nu) / 2}
    +
    \frac{2 \sigma}{D} \Bigl( \frac{1}{3} k (k + 1)^2 \Bigr)^{1 / 2}
    \\
    &=
    \frac{
      2^{(1 + 3 \nu) / 2} (1 + \nu)^{(1 + \nu) / 2}
    }{
      3^{(1 - \nu) / 2}
    }
    \frac{L_{\nu}}{D^{1 - \nu}}
    k^{(1 - \nu) / 2} (k + 1)^{1 - \nu}
    +
    \frac{2}{\sqrt{3}} \frac{\sigma}{D} \sqrt{k} \, (k + 1)
    \\
    &\leq
    \frac{8 L_{\nu}}{D^{1 - \nu}}
    k^{(1 - \nu) / 2} (k + 1)^{1 - \nu}
    +
    \frac{2}{\sqrt{3}} \frac{\sigma}{D} \sqrt{k} \, (k + 1),
  \end{align*}
  where we have used the fact that
  $
    \sum_{i = 1}^k i^2
    =
    \frac{1}{6} k (k + 1) (2 k + 1)
    \leq
    \frac{1}{3} k (k + 1)^2
  $.
  Consequently,
  \begin{align*}
    \frac{4 \Expectation [H_k] D^2}{k (k + 1)}
    &\leq
    32 L_{\nu} D^{1 + \nu}
    \frac{k^{(1 - \nu) / 2} (k + 1)^{1 - \nu}}{k (k + 1)}
    +
    \frac{8 \sigma D}{\sqrt{3}}
    \frac{\sqrt{k} \, (k + 1)}{k (k + 1)}
    \\
    &=
    \frac{32 L_{\nu} D^{1 + \nu}}{k^{(1 + \nu) / 2} (k + 1)^{\nu}}
    +
    \frac{8 \sigma D}{\sqrt{3 k}}
    \leq
    \frac{32 L_{\nu} D^{1 + \nu}}{k^{(1 + 3 \nu) / 2}}
    +
    \frac{8 \sigma D}{\sqrt{3 k}}.
  \end{align*}
  This proves the final inequality in \cref{\LocalName{eq:Result}}
  since $\nu \in \ClosedClosedInterval{0}{1}$ was arbitrary.
\end{proof}

%% file: Appendix/FastGradientMethod.tex
\section{Universal Line-Search-Free Fast Gradient Method}

\begin{theorem}
  \label{th:ConvergenceRateOfFastGradientMethod}
  \UsingNamespace{ConvergenceRateOfFastGradientMethod}
  Let \cref{alg:StochasticFastGradientMethod} be applied for solving
  problem~\eqref{eq:MainProblem} under
  \cref{as:BoundedFeasibleSet,as:AtLeastOneHolderConstantIsFinite} with the
  deterministic gradient oracle $g(x, \xi) \equiv \nabla f(x)$ and with
  $
    \hat{\beta}_{k + 1}
    \DefinedEqual
    \BregmanDistanceWithSubgradient{f}{g_k^y}(y_k, x_{k + 1})
  $
  (\cref{StochasticFastGradientMethod::step:UpdateRegularizationParameter})
  at each iteration $k \geq 0$.
  Then, for all $k \geq 1$, it holds that
  \begin{equation}
    \LocalLabel{eq:Result}
    F(x_k) - F^*
    \leq
    \frac{4 H_k D^2}{k (k + 1)}
    \leq
    \inf_{\nu \in \ClosedClosedInterval{0}{1}}
    \frac{8 L_{\nu} D^{1 + \nu}}{k^{(1 + 3 \nu) / 2}}.
  \end{equation}
\end{theorem}

\begin{proof}
  \UsingNestedNamespace{ConvergenceRateOfFastGradientMethod}{Proof}
  We proceed exactly in the same way as in the proof of
  \cref{th:ConvergenceRateOfStochasticFastGradientMethod}
  (\cref{sec:ProofForStochasticFastGradientMethod})
  but do not upper bound $\beta_{k + 1} = \BregmanDistance{f}(y_k, x_{k + 1})$
  with $\hat{\beta}_{k + 1}$ in
  \cref{ConvergenceRateOfStochasticFastGradientMethod::Proof::eq:UpperBoundOnBregmanDistanceApproximation}.
  We then arrive, exactly as before, at the following inequality
  that holds for any $k \geq 1$:
  \begin{equation}
    \LocalLabel{eq:BoundForFunctionResidual}
    F(x_k) - F^* \leq \frac{4 H_k D^2}{k (k + 1)}.
  \end{equation}

  To upper bound~$H_k$, we use, as before, the following equation:
  \[
    (H_{k + 1} - H_k) D^2
    =
    \PositivePart{
      A_{k + 1} \beta_{k + 1} - \tfrac{1}{2} H_{k + 1} r_{k + 1}^2
    },
  \]
  that holds for any $k \geq 0$, but now we can upper bound
  \[
    \beta_{k + 1}
    \leq
    \frac{L_{\nu}}{1 + \nu} A_{k + 1} \Norm{x_{k + 1} - y_k}^{1 + \nu}.
  \]
  This is essentially the same bound that we had in
  \cref{ConvergenceRateOfStochasticFastGradientMethod::Proof::eq:UpperBoundOnStochasticSymmetrizedBregmanDistance}
  with the formal change of $L_{\nu}$ to
  $L_{\nu}' \DefinedEqual \frac{L_{\nu}}{1 + \nu}$.
  Proceeding exactly as before, we then obtain
  \[
      A_{k + 1} \beta_{k + 1}
      \leq
      \frac{2^{\nu} L_{\nu}}{1 + \nu} (k + 1)^{1 - \nu} r_{k + 1}^{1 + \nu}
  \]
  for any $k \geq 0$ and arbitrary $\nu \in \ClosedClosedInterval{0}{1}$,
  which gives us
  \[
      (H_{k + 1} - H_k) D^2
      \leq
      \PositivePart[\Big]{
          \frac{2^{\nu} L_{\nu}}{1 + \nu} (k + 1)^{1 - \nu} r_{k + 1}^{1 + \nu}
          -
          \frac{H_{k + 1}}{2} r_{k + 1}^2
      }.
  \]

  Instead of
  \cref{th:EstimatingGrowthRateOfRegularizationParameterInStochasticCase},
  we can now apply a slightly more precise result (in terms of absolute
  constants)---\cref{th:EstimatingGrowthRateOfRegularizationParameter}
  (with $\Omega \DefinedEqual D^2$, $M \DefinedEqual 2^{\nu} L_{\nu}$,
  $\gamma_k \DefinedEqual k$)---to conclude that, for all $k \geq 1$,
  \begin{align*}
    H_k
    \leq
    \Bigl[
      \frac{1}{(1 + \nu) D^2} \sum_{i = 1}^k i^2
    \Bigr]^{(1 - \nu) / 2}
    2^{\nu} L_{\nu}
    &\leq
    2^{\nu}
    \Bigl[ \frac{1}{3 (1 + \nu)} k (k + 1)^2 \Bigr]^{(1 - \nu) / 2}
    \frac{L_{\nu}}{D^{1 - \nu}}
    \\
    &\leq
    \frac{2 L_{\nu}}{D^{1 - \nu}}
    k^{(1 - \nu) / 2} (k + 1)^{1 - \nu},
  \end{align*}
  where the second inequality is due to
  $
    \sum_{i = 1}^k i^2
    =
    \frac{1}{6} k (k + 1) (2 k + 1)
    \leq
    \frac{1}{3} k (k + 1)^2
  $,
  and the final inequality follows from the fact that
  $2^{\nu} / [3 (1 + \nu)]^{(1 - \nu) / 2}$ monotonically
  increases in $\nu \in \ClosedClosedInterval{0}{1}$.
  Substituting the above bound into
  \cref{\LocalName{eq:BoundForFunctionResidual}}, we get
  \[
    F(x_k) - F^*
    \leq
    2 L_{\nu} D^{1 + \nu}
    \frac{k^{(1 - \nu) / 2} (k + 1)^{1 - \nu}}{k (k + 1)}
    =
    \frac{
      2 L_{\nu} D^{1 + \nu}
    }{
      k^{(1 + \nu) / 2} (k + 1)^{\nu}
    }
    \leq
    \frac{2 L_{\nu} D^{1 + \nu}}{k^{(1 + 3 \nu) / 2}}.
    \qedhere
  \]
\end{proof}

%% file: Appendix/AuxiliaryResults.tex
\section{Auxiliary Results}

\begin{lemma}
  \label{th:SolvingEquationForRegularizationParameter}
  \UsingNamespace{SolvingEquationForRegularizationParameter}
  Let $H, \beta, \rho \geq 0$ and $\Omega > 0$.
  Then, the equation
  \begin{equation}
    \LocalLabel{eq:Equation}
    (H_+ - H) \Omega = \PositivePart{\beta - H_+ \rho}
  \end{equation}
  has a unique solution given by
  \begin{equation}
    \LocalLabel{eq:Solution}
    H_+
    \DefinedEqual
    H
    +
    \frac{\PositivePart{\beta - H \rho}}{\Omega + \rho}.
  \end{equation}
\end{lemma}

\begin{proof}
  \UsingNestedNamespace{SolvingEquationForRegularizationParameter}{Proof}

  Denote the left- and right-hand sides in \cref{\MasterName{eq:Equation}} (as
  functions of $H_+$) by $\zeta_1(H_+)$ and $\zeta_2(H_+)$, respectively, and
  let $\zeta(H_+) \DefinedEqual \zeta_1(H_+) - \zeta_2(H_+)$.
  Note that both $\zeta_1$ and $\zeta_2$ are continuous functions, $\zeta_1$ is
  strictly increasing, while $\zeta_2$ is decreasing, hence $\zeta$ is a
  continuous strictly increasing function.
  When $H_+ = H$, we have $\zeta_1(H) = 0$, while $\zeta_2(H) \geq 0$, hence
  $\zeta(H) \leq 0$.
  When $H_+ \to +\infty$, $\zeta_1(H_+)$ tends to $+\infty$, while
  $\zeta_2(H_+)$ tends to a finite number (either $0$ if $\rho > 0$, or $\beta$
  if $\rho = 0$), hence $\zeta(H_+)$ tends to $+\infty$.
  Thus, there exists a unique point $H_+ \geq H$ such that $\zeta(H_+) = 0$.
  This point is exactly the unique solution of
  equation~\eqref{\MasterName{eq:Equation}}.

  It remains to show that~\eqref{\MasterName{eq:Solution}} is indeed a solution
  to \cref{\MasterName{eq:Equation}}.
  But this is simple.
  Indeed, if $\beta \leq H \rho$, then, according to
  \cref{\MasterName{eq:Solution}}, $H_+ = H + (\beta - H \rho) / (\Omega + \rho)$,
  and hence
  \[
    \beta - H_+ \rho
    =
    \beta
    -
    H \rho
    -
    \frac{\beta - H \rho}{\Omega + \rho} \rho
    =
    \frac{\Omega}{\Omega + \rho} (\beta - H \rho)
    =
    (H_+ - H) \Omega
    \quad (\geq 0),
  \]
  which means that $H_+$ satisfies \cref{\MasterName{eq:Equation}}.
  If $\beta > H \rho$, then, by \cref{\MasterName{eq:Solution}}, $H_+ = H$, and
  hence
  \[
    \PositivePart{\beta - H_+ \rho}
    =
    \PositivePart{\beta - H \rho}
    =
    0
    =
    (H_+ - H) \Omega,
  \]
  which also means that $H_+$ satisfies \cref{\MasterName{eq:Equation}}.
\end{proof}

\begin{lemma}
  \label{th:OptimalityConditionForProximalPointStep}

  Let $\Map{\zeta}{\RealField^n}{\RealFieldPlusInfty}$ be a convex function,
  $\bar{x} \in \EffectiveDomain \zeta$, $H \geq 0$.
  Then, for any $x^* \in \ArgminSet_{x \in \EffectiveDomain \zeta} \{ \zeta(x) +
  \tfrac{1}{2} H \Norm{x -\bar{x}}^2 \}$ and any $x \in \EffectiveDomain \zeta$,
  we have
  \[
    \zeta(x) + \tfrac{1}{2} H \Norm{x - \bar{x}}^2
    \geq
    \zeta(x^*) + \tfrac{1}{2} H \Norm{x^* - \bar{x}}^2
    +
    \tfrac{1}{2} H \Norm{x - x^*}^2.
  \]
\end{lemma}

\begin{proof}
  This is a standard result that can be seen as a consequence of the fact that
  $\zeta_H(x) \DefinedEqual \zeta(x) + \frac{H}{2} \Norm{x - \bar{x}}^2$ is a
  strongly convex function with constant~$H$, and hence
  $\zeta_H(x) \geq \zeta_H(x^*) + \frac{H}{2} \Norm{x - x^*}^2$
  for any $x \in \EffectiveDomain \zeta$.
\end{proof}

\begin{lemma}
  \label{th:UpperBoundOnRegularizedBregmanDistanceTerm}
  \UsingNamespace{UpperBoundOnRegularizedBregmanDistanceTerm}
  Let $\nu \in \ClosedOpenInterval{0}{1}$, $M \geq 0$ and $H > 0$.
  Then,
  \begin{equation}
    \LocalLabel{eq:Result}
    \max_{r \geq 0} \Bigl\{
      \frac{M}{1 + \nu} r^{1 + \nu}
      -
      \frac{H}{2} r^2
    \Bigr\}
    =
    \frac{1 - \nu}{2 (1 + \nu)}
    \frac{M^{2 / (1 - \nu)}}{H^{(1 + \nu) / (1 - \nu)}}.
  \end{equation}
\end{lemma}

\begin{proof}
  \UsingNestedNamespace{UpperBoundOnRegularizedBregmanDistanceTerm}{Proof}

  After the change of variables $t = r^{1 + \nu}$, the objective function inside
  the $\max$ becomes concave in~$t$
  (since $r^2 = t^{2 / (1 + \nu)}$ with $\frac{2}{1 + \nu} \geq 1$).
  Computing its derivative and setting to zero, we see that the maximum is
  attained at the point $r_* \DefinedEqual (M / H)^{1  / (1 - \nu)}$.
  Thus,
  \begin{align*}
    \hspace{2em}&\hspace{-2em}
    \max_{r \geq 0} \Bigl\{
      \frac{M}{1 + \nu} r^{1 + \nu}
      -
      \frac{H}{2} r^2
    \Bigr\}
    =
    \frac{M}{1 + \nu} \Bigl( \frac{M}{H} \Bigr)^{(1 + \nu) / (1 - \nu)}
    -
    \frac{H}{2} \Bigl( \frac{M}{H} \Bigr)^{2 / (1 - \nu)}
    \\
    &=
    \frac{1}{1 + \nu}
    \frac{M^{2 / (1 - \nu)}}{H^{(1 + \nu) / (1 - \nu)}}
    \Bigl( 1 - \frac{1}{2} (1 + \nu) \Bigr)
    =
    \frac{1 - \nu}{2 (1 + \nu)}
    \frac{M^{2 / (1 - \nu)}}{H^{(1 + \nu) / (1 - \nu)}}.
    \qedhere
  \end{align*}
\end{proof}

\begin{lemma}
  \label{th:RecurrentInequalityWithPowerGrowth}
  Let $(H_k)_{k = 0}^\infty$ be a nonnegative nondecreasing sequence of reals
  such that, for any $k \geq 0$,
  \[
    (p + 1) H_{k + 1}^p (H_{k + 1} - H_k) \leq \alpha_{k + 1},
  \]
  where $p \geq 0$ is real and $(\alpha_k)_{k = 1}^\infty$ is a nonnegative
  sequence of reals.
  Then, for any $k \geq 1$, it holds that
  \[
    H_k \leq \Bigl( H_0^{p + 1} + \sum_{i = 1}^k \alpha_i \Bigr)^{1 / (p + 1)}.
  \]
\end{lemma}

\begin{proof}
  Since $H_k \leq H_{k + 1}$ for any $k \geq 0$ and $p \geq 0$, we can estimate
  \[
    \alpha_{k + 1}
    \geq
    (p + 1) H_{k + 1}^p (H_{k + 1} - H_k)
    \geq
    (p + 1) \int_{H_k}^{H_{k + 1}} t^p dt
    =
    H_{k + 1}^{p + 1} - H_k^{p + 1}.
  \]
  Telescoping these inequalities, we obtain, for any $k \geq 1$,
  \[
    H_k^{p + 1} - H_0^{p + 1} \leq \sum_{i = 1}^k \alpha_i,
  \]
  and the claim follows.
\end{proof}

\begin{lemma}
  \label{th:RecurrentInequalityWithQuadraticGrowth}
  Let $(H_k)_{k = 0}^\infty$ be a nonnegative nondecreasing sequence of reals
  such that, for any $k \geq 0$,
  \[
    H_{k + 1} (H_{k + 1} - H_k) \geq \alpha_{k + 1},
  \]
  where $(\alpha)_{k = 1}^\infty$ is a nonnegative sequence of reals.
  Then, for any $k \geq 0$, it holds that
  \[
    H_k \geq \Bigl( H_0^2 + \sum_{i = 1}^k \alpha_i \Bigr)^{1 / 2}.
  \]
\end{lemma}

\begin{proof}
  Indeed, for any $k \geq 0$, we can estimate
  \[
    \alpha_{k + 1}
    \leq
    H_{k + 1} (H_{k + 1} - H_k)
    \leq
    (H_{k + 1} + H_k) (H_{k + 1} - H_k)
    =
    H_{k + 1}^2 - H_k^2.
  \]
  Summing up these inequalities and rearranging, we obtain the claim.
\end{proof}

\begin{lemma}
  \label{th:EstimatingGrowthRateOfRegularizationParameter}
  \UsingNamespace{EstimatingGrowthRateOfRegularizationParameter}
  Let $(H_k)_{k = 0}^\infty$ be a nondecreasing sequence
  such that $H_0 = 0$ and, for all $k \geq 0$, it holds
  \begin{equation}
    \LocalLabel{eq:RecurrenceRelation}
    (H_{k + 1} - H_k) \Omega
    \leq
    \PositivePart[\Big]{
      \frac{1}{1 + \nu} M \gamma_{k + 1}^{1 - \nu} r_{k + 1}^{1 + \nu}
      -
      \frac{1}{2} H_{k + 1} r_{k + 1}^2
    },
  \end{equation}
  where $\Omega > 0$, $M \geq 0$, $\nu \in \ClosedClosedInterval{0}{1}$ are
  certain constants, and $(\gamma_k)_{k = 1}^\infty$ and
  $(r_k)_{k = 1}^\infty$ are certain positive and nonnegative
  sequences, respectively.
  Then, for all $k \geq 1$, we have
  \begin{equation}
    \LocalLabel{eq:ExplicitRate}
    H_k
    \leq
    \Bigl[
      \frac{1}{(1 + \nu) \Omega} \sum_{i = 1}^k \gamma_i^2
    \Bigr]^{(1 - \nu) / 2}
    M.
  \end{equation}
\end{lemma}

\begin{proof}
  \UsingNestedNamespace{EstimatingGrowthRateOfRegularizationParameter}{Proof}

  Suppose $\nu = 1$.
  Then, according to \cref{\MasterName{eq:RecurrenceRelation}},
  for all $k \geq 0$, we have
  \begin{equation}
    \LocalLabel{eq:LipschitzCase:RecurrenceRelation}
    (H_{k + 1} - H_k) \Omega
    \leq
    \PositivePart{
      \tfrac{1}{2} M r_{k + 1}^2 - \tfrac{1}{2} H_{k + 1} r_{k + 1}^2
    }
    =
    \PositivePart{M - H_{k + 1}} \tfrac{1}{2} r_{k + 1}^2.
  \end{equation}
  Since $H_0 = 0 \leq M$, this implies that $H_k \leq M$ for all $k \geq 0$
  (which is exactly \cref{\MasterName{eq:ExplicitRate}} for $\nu = 1$).
  Indeed, if $H_k \leq M < H_{k + 1}$ for some $k \geq 0$, then the left-hand
  side in \cref{\LocalName{eq:LipschitzCase:RecurrenceRelation}} is strictly
  positive, while the right-hand side is zero, which is a contradiction.

  From now on, suppose $\nu < 1$.
  Without loss of generality, we can assume that $H_{k + 1} > 0$
  for all $k \geq 0$.
  Let $k \geq 0$ be arbitrary.
  Applying \cref{th:UpperBoundOnRegularizedBregmanDistanceTerm} to bound the
  right-hand side in \cref{\MasterName{eq:RecurrenceRelation}} (and using the
  monotonicity of $\PositivePart{\cdot}$), we obtain
  \[
    (H_{k + 1} - H_k) \Omega
    \leq
    \frac{1 - \nu}{2 (1 + \nu)}
    \frac{
      (M \gamma_{k + 1}^{1 - \nu})^{2 / (1 - \nu)}
    }{
      H_{k + 1}^{(1 + \nu) / (1 - \nu)}
    }
    =
    \frac{1 - \nu}{2 (1 + \nu)}
    \frac{
      M^{2 / (1 - \nu)}
    }{
      H_{k + 1}^{(1 + \nu) / (1 - \nu)}
    }
    \gamma_{k + 1}^2.
  \]
  Applying \cref{th:RecurrentInequalityWithPowerGrowth}
  (with $p = \frac{1 + \nu}{1 - \nu}$ for which $p + 1 = \frac{2}{1 - \nu}$)
  and using the fact that $H_0 = 0$, we conclude that
  \[
    H_k
    \leq
    \Bigl[
      \frac{M^{2 / (1 - \nu)}}{2 (1 + \nu) \Omega} \sum_{i = 1}^k \gamma_i^2
    \Bigr]^{(1 - \nu) / 2}
    =
    \Bigl[ \frac{1}{(1 + \nu) \Omega} \sum_{i = 1}^k \gamma_i^2 \Bigr]^{(1 - \nu) / 2} M.
    \qedhere
  \]
\end{proof}

\begin{lemma}
  \label{th:EstimatingGrowthRateOfRegularizationParameterInStochasticCase}
  \UsingNamespace{EstimatingGrowthRateOfRegularizationParameterInStochasticCase}
  Let $(H_k)_{k = 0}^\infty$ be a nondecreasing sequence
  such that $H_0 = 0$ and, for all $k \geq 0$, it holds
  \begin{equation}
    \LocalLabel{eq:RecurrenceRelation}
    (H_{k + 1} - H_k) \Omega
    \leq
    \PositivePart{
      L \alpha_{k + 1}^{1 - \nu} r_{k + 1}^{1 + \nu}
      +
      \gamma_{k + 1} r_{k + 1}
      -
      \tfrac{1}{2} H_{k + 1} r_{k + 1}^2
    },
  \end{equation}
  where $\Omega > 0$, $M \geq 0$, $\nu \in \ClosedClosedInterval{0}{1}$ are
  certain constants, $(\alpha_k)_{k = 1}^{\infty}$ is a certain positive
  sequence, and $(r_k)_{k = 1}^\infty$ and $(\gamma_k)_{k = 1}^\infty$ are
  certain nonnegative sequences.
  Then, for all $k \geq 1$,
  \begin{equation}
    \LocalLabel{eq:ExplicitRate}
    H_k
    \leq
    [2 (1 + \nu)]^{(1 + \nu) / 2}
    L
    \Bigl( \frac{1}{\Omega} \sum_{i = 1}^k \alpha_i^2 \Bigr)^{(1 - \nu) / 2}
    +
    \Bigl( \frac{2}{\Omega} \sum_{i = 1}^k \gamma_i^2 \Bigr)^{1 / 2}.
  \end{equation}
\end{lemma}

\begin{remark}
  Setting $\gamma_k \equiv 0$ in
  \cref{th:EstimatingGrowthRateOfRegularizationParameterInStochasticCase},
  we recover
  \cref{th:EstimatingGrowthRateOfRegularizationParameter}.
\end{remark}

\begin{proof}
  \UsingNestedNamespace{EstimatingGrowthRateOfRegularizationParameterInStochasticCase}{Proof}

  \ProofPart

  Without loss of generality, we can assume that $H_{k + 1} > 0$
  for all $k \geq 0$.
  Indeed, otherwise, either $H_k = 0$ for all $k \geq 0$, and
  \cref{\MasterName{eq:ExplicitRate}} is trivial, or we can work with the
  subsequence $(H_k)_{k = k_0}^\infty$, where $k_0 \geq 0$ is the first integer
  such that $H_{k_0 + 1} > 0$.

  \ProofPart

  Suppose%
  \footnote{
    In principle, we can cover the case $\nu = 1$ by only considering the values
    of $\nu \in \ClosedOpenInterval{0}{1}$ and then passing to the limit as $\nu
    \to 1$.
    However, we prefer to present a more explicit proof without using the
    limiting argument.
  }
  $\nu = 1$.
  In this case, \cref{\MasterName{eq:RecurrenceRelation}} reads
  \begin{equation}
    \LocalLabel{eq:LipschitzCase:RecurrenceRelation}
    (H_{k + 1} - H_k) \Omega
    \leq
    \PositivePart{
      (L - \tfrac{1}{2} H_{k + 1}) r_{k + 1}^2
      +
      \gamma_{k + 1} r_{k + 1}
    }
  \end{equation}
  for all $k \geq 0$, and we need to prove that, for all $k \geq 1$,
  \begin{equation}
    \LocalLabel{eq:LipschitzCase:ExplicitRate}
    H_k
    \leq
    4 L
    +
    \Bigl( \frac{2}{\Omega} \sum_{i = 1}^k \gamma_i^2 \Bigr)^{1 / 2}.
  \end{equation}
  (This is exactly \cref{\MasterName{eq:ExplicitRate}} for $\nu = 1$.)

  Since $H_0 = 0$, we can assume that there exists an index $k_0 \geq 0$
  such that
  \begin{equation}
    \LocalLabel{eq:LipschitzCase:SwitchingMoment}
    H_{k_0} \leq 4 L < H_{k_0 + 1}.
  \end{equation}
  (Otherwise, $H_k \leq 4 L$ for all $k \geq 0$, and
  \cref{\LocalName{eq:LipschitzCase:ExplicitRate}} is trivial.)
  As $(H_k)_{k = 0}^\infty$ is nondecreasing,
  \cref{\LocalName{eq:LipschitzCase:ExplicitRate}} is clearly valid for all
  indices $0 \leq k \leq k_0$.
  Let us prove that it is also valid for all $k \geq k_0 + 1$.

  Let $k \geq k_0$ be arbitrary.
  By monotonicity of $(H_i)_{i = 0}^\infty$, from
  \cref{\LocalName{eq:LipschitzCase:SwitchingMoment}}, it follows that
  $H_{k + 1} \geq H_{k_0 + 1} > 4 L$.
  Therefore,
  \[
    \Bigl( L - \frac{1}{2} H_{k + 1} \Bigr) r_{k + 1}^2
    +
    \gamma_{k + 1} r_{k + 1}
    \leq
    \gamma_{k + 1} r_{k + 1}
    -
    \frac{1}{4} H_{k + 1} r_{k + 1}^2
    \leq
    \frac{\gamma_{k + 1}^2}{H_{k + 1}},
  \]
  where the final inequality follows from
  \cref{th:UpperBoundOnRegularizedBregmanDistanceTerm}
  (with $\nu \DefinedEqual 0$ and $H \DefinedEqual \frac{1}{2} H_{k + 1}$).
  Combining this with \cref{\LocalName{eq:LipschitzCase:RecurrenceRelation}}
  (using the monotonicity of $\PositivePart{\cdot}$), we get
  \[
    (H_{k + 1} - H_k) \Omega
    \leq
    \frac{\gamma_{k + 1}^2}{H_{k + 1}}.
  \]
  Thus, for all $k \geq k_0$, we have
  \begin{equation}
    \LocalLabel{eq:LipschitzCase:SimplifiedRecurrenceRelation}
    (H_{k + 1}^2 - H_k^2) \Omega
    \leq
    2 H_{k + 1} (H_{k + 1} - H_k) \Omega
    \leq
    2 \gamma_{k + 1}^2.
  \end{equation}
  (Recall that $H_k \leq H_{k + 1}$.)

  Let $k \geq k_0 + 1$ be arbitrary.
  Summing up \cref{\LocalName{eq:LipschitzCase:SimplifiedRecurrenceRelation}}
  for all indices $k_0 \leq k' \leq k - 1$ and rearranging, we get
  \[
    H_k^2
    \leq
    H_{k_0}^2
    +
    \frac{2}{\Omega} \sum_{i = k_0 + 1}^k \gamma_i^2
    \leq
    (4 L)^2
    +
    \frac{2}{\Omega} \sum_{i = 1}^k \gamma_i^2,
  \]
  where the last inequality is due to
  \cref{\LocalName{eq:LipschitzCase:SwitchingMoment}}
  (and the fact that $k_0 \geq 0$).
  Using the fact that $\sqrt{a + b} \leq \sqrt{a} + \sqrt{b}$
  for any $a, b \geq 0$, we obtain
  \cref{\LocalName{eq:LipschitzCase:ExplicitRate}}.

  \ProofPart

  Now suppose $\nu < 1$.
  Let $k \geq 0$ be arbitrary.
  Applying \cref{th:UpperBoundOnRegularizedBregmanDistanceTerm} twice,
  we obtain
  \begin{align*}
    \hspace{2em}&\hspace{-2em}
    L \alpha_{k + 1}^{1 - \nu} r_{k + 1}^{1 + \nu}
    +
    \gamma_{k + 1} r_{k + 1}
    -
    \tfrac{1}{2} H_{k + 1} r_{k + 1}^2
    \\
    &=
    [
      L \alpha_{k + 1}^{1 - \nu} r_{k + 1}^{1 + \nu}
      -
      \tfrac{1}{4} H_{k + 1} r_{k + 1}^2
    ]
    +
    [\gamma_{k + 1} r_{k + 1} - \tfrac{1}{4} H_{k + 1} r_{k + 1}^2]
    \\
    &\leq
    \frac{1 - \nu}{2 (1 + \nu)}
    \frac{
      [(1 + \nu) L \alpha_{k + 1}^{1 - \nu}]^{2 / (1 - \nu)}
    }{
      ( \frac{1}{2} H_{k + 1} )^{(1 + \nu) / (1 - \nu)}
    }
    +
    \frac{1}{2} \frac{\gamma_{k + 1}^2}{\frac{1}{2} H_{k + 1}}
    \\
    &=
    (1 - \nu)
    \frac{M^{2 / (1 - \nu)}}{H_{k + 1}^{(1 + \nu) / (1 - \nu)}}
    \alpha_{k + 1}^2
    +
    \frac{\gamma_{k + 1}^2}{H_{k + 1}},
  \end{align*}
  where
  \begin{equation}
    \LocalLabel{eq:ScaledHolderConstant}
    M
    \DefinedEqual
    \frac{
      2^{(1 + \nu) / 2} (1 + \nu)
    }{
      [2 (1 + \nu)]^{(1 - \nu) / 2}
    }
    L
    =
    2^{\nu} (1 + \nu)^{(1 + \nu) / 2} L.
  \end{equation}
  Combining this with \cref{\MasterName{eq:RecurrenceRelation}}
  (using the monotonicity of $\PositivePart{\cdot}$), we get
  \[
    (H_{k + 1} - H_k) \Omega
    \leq
    (1 - \nu)
    \frac{M^{2 / (1 - \nu)}}{H_{k + 1}^{(1 + \nu) / (1 - \nu)}}
    \alpha_{k + 1}^2
    +
    \frac{\gamma_{k + 1}^2}{H_{k + 1}}.
  \]
  Since $H_k \leq H_{k + 1}$, it follows that
  \[
    \frac{1}{2} (H_{k + 1}^2 - H_k^2) \Omega
    \leq
    H_{k + 1} (H_{k + 1} - H_k) \Omega
    \leq
    (1 - \nu)
    \frac{M^{2 / (1 - \nu)}}{H_{k + 1}^{2 \nu / (1 - \nu)}}
    \alpha_{k + 1}^2
    +
    \gamma_{k + 1}^2.
  \]
  Note that this inequality is valid for all $k \geq 0$.

  Applying \cref{th:EstimatingGrowthRateForSimplifiedRecurrenceInStochasticCase}
  (with
  $C_k \DefinedEqual H_k^2$,
  $\alpha_k' \DefinedEqual \frac{2}{\Omega} \alpha_k^2$
  and~$\gamma_k' \DefinedEqual \frac{2}{\Omega} \gamma_k^2$),
  we conclude that, for all $k \geq 1$,
  \begin{align*}
    H_k^2
    &\leq
    M^2 \Bigl( \sum_{i = 1}^k \frac{2}{\Omega} \alpha_i^2 \Bigr)^{1 - \nu}
    +
    \sum_{i = 1}^k \frac{2}{\Omega} \gamma_i^2
    =
    2^{2 \nu} (1 + \nu)^{1 + \nu} L^2
    \Bigl( \sum_{i = 1}^k \frac{2}{\Omega} \alpha_i^2 \Bigr)^{1 - \nu}
    +
    \frac{2}{\Omega} \sum_{i = 1}^k \gamma_i^2
    \\
    &=
    [2 (1 + \nu)]^{1 + \nu} L^2
    \Bigl( \frac{1}{\Omega} \sum_{i = 1}^k \alpha_i^2 \Bigr)^{1 - \nu}
    +
    \frac{2}{\Omega} \sum_{i = 1}^k \gamma_i^2,
  \end{align*}
  where the second identity follows from
  \cref{\LocalName{eq:ScaledHolderConstant}}.
  Using the fact that $\sqrt{a + b} \leq \sqrt{a} + \sqrt{b}$
  for any $a, b \geq 0$, we obtain \cref{\MasterName{eq:ExplicitRate}}.
\end{proof}

\begin{lemma}
  \label{th:EstimatingGrowthRateForSimplifiedRecurrenceInStochasticCase}
  \UsingNamespace{EstimatingGrowthRateForSimplifiedRecurrenceInStochasticCase}
  Let $(C_k)_{k = 1}^\infty$ be a positive sequence satisfying,
  for all $k \geq 0$,
  \begin{equation}
    \LocalLabel{eq:RecurrenceRelation}
    C_{k + 1} - C_k
    \leq
    (1 - \nu)
    \frac{M^{2 / (1 - \nu)}}{C_{k + 1}^{\nu / (1 - \nu)}}
    \alpha_{k + 1}
    +
    \gamma_{k + 1},
  \end{equation}
  where $C_0 \DefinedEqual 0$,
  and $M \geq 0$, $\nu \in \ClosedOpenInterval{0}{1}$ are certain constants,
  and $(\alpha_k)_{k = 1}^\infty$ and~$(\gamma_k)_{k = 1}^{\infty}$
  are certain nonnegative sequences.
  Then, for all $k \geq 1$, we have
  \begin{equation}
    \LocalLabel{eq:ExplicitRate}
    C_k
    \leq
    M^2
    \Bigl( \sum_{i = 1}^k \alpha_i \Bigr)^{1 - \nu}
    +
    \sum_{i = 1}^k \gamma_i.
  \end{equation}
\end{lemma}

\begin{proof}
  \UsingNestedNamespace{EstimatingGrowthRateForSimplifiedRecurrenceInStochasticCase}{Proof}

  For each $k \geq 0$, let $\hat{C}_k$ be the right-hand side of
  \cref{\MasterName{eq:ExplicitRate}}:
  \begin{equation}
    \LocalLabel{eq:UpperBoundOnRate}
    \hat{C}_k
    \DefinedEqual
    M^2 A_k^{1 - \nu} + \sum_{i = 1}^k \gamma_i,
    \qquad
    A_k \DefinedEqual \sum_{i = 1}^k \alpha_i,
  \end{equation}
  with the convention that $\hat{C}_0 = A_0 = 0$.
  Note that $\hat{C}_k > 0$ for all $k \geq 1$.
  Indeed, if $\hat{C}_k = 0$ for some $k \geq 1$, then $\hat{C}_1 = 0$
  (by the monotonicity of~$(\hat{C}_k)_{k = 1}^{\infty}$),
  which means that $\gamma_1 = 0$ and either $M = 0$ or $\alpha_1 = 0$;
  but then, according to \cref{\MasterName{eq:RecurrenceRelation}},
  $C_1 - C_0 \leq 0$;
  since $C_0 = 0$, this implies $C_1 \leq 0$, which contradicts our assumption
  about the positivity of~$(C_k)_{k = 1}^\infty$.

  Let us prove by induction that $C_k \leq \hat{C}_k$ for all $k \geq 0$.
  Clearly, this inequality is satisfied for $k = 0$ since $C_0 = \hat{C}_0 = 0$.
  Now suppose that $C_k \leq \hat{C}_k$ for some $k \geq 0$, and let us prove
  that $C_{k + 1} \leq \hat{C}_{k + 1}$.

  Let
  $\Map{\chi_{k + 1}}{\OpenOpenInterval{0}{+\infty}}{\RealField}$
  be the function
  \begin{equation}
    \LocalLabel{eq:ResolventFunction}
    \chi_{k + 1}(C)
    \DefinedEqual
    C
    -
    (1 - \nu) \frac{M^{2 / (1 - \nu)}}{C^{\nu / (1 - \nu)}} \alpha_{k + 1}.
  \end{equation}
  According to \cref{\MasterName{eq:RecurrenceRelation}} and the inductive
  hypothesis, we have
  \begin{equation}
    \LocalLabel{eq:PreliminaryEstimate}
    \chi_{k + 1}(C_{k + 1})
    \leq
    C_k + \gamma_{k + 1}
    \leq
    \hat{C}_k + \gamma_{k + 1}.
  \end{equation}
  Since the function~$\chi_{k + 1}$ is strictly increasing, to prove that
  $C_{k  + 1} \leq \hat{C}_{k + 1}$, it suffices to show that
  $\chi_{k + 1}(C_{k + 1}) \leq \chi_{k + 1}(\hat{C}_{k + 1})$.
  According to \cref{\LocalName{eq:PreliminaryEstimate}}, for this, it suffices
  to show that
  \[
    \hat{C}_k + \gamma_{k + 1}
    \leq
    \chi_{k + 1}(\hat{C}_{k + 1}).
  \]
  Substituting \cref{\LocalName{eq:ResolventFunction}} and rearranging, we
  see that we need to prove that $(\hat{C}_i)_{i = 0}^\infty$
  satisfies \cref{\MasterName{eq:ExplicitRate}} with the reversed sign:
  \[
    \hat{C}_{k + 1} - \hat{C}_k
    \geq
    (1 - \nu)
    \frac{M^{2 / (1 - \nu)}}{\hat{C}_{k + 1}^{\nu / (1 - \nu)}}
    \alpha_{k + 1}
    +
    \gamma_{k + 1}.
  \]

  In view of \cref{\LocalName{eq:UpperBoundOnRate}}, we have
  \[
    \hat{C}_{k + 1} - \hat{C}_k
    =
    M^2 [A_{k + 1}^{1 - \nu} - A_k^{1 - \nu}]
    +
    \gamma_{k + 1}.
  \]
  Thus, we need to check if
  \[
    M^2 [A_{k + 1}^{1 - \nu} - A_k^{1 - \nu}]
    \geq
    (1 - \nu)
    \frac{M^{2 / (1 - \nu)}}{\hat{C}_{k + 1}^{\nu / (1 - \nu)}}
    \alpha_{k + 1},
  \]
  or, equivalently, if
  \[
    \hat{C}_{k + 1}^{\nu / (1 - \nu)}
    [A_{k + 1}^{1 - \nu} - A_k^{1 - \nu}]
    \geq
    (1 - \nu) M^{2 \nu / (1 - \nu)} \alpha_{k + 1}.
  \]
  From \cref{\LocalName{eq:UpperBoundOnRate}}, it follows that
  $\hat{C}_{k + 1} \geq M^2 A_{k + 1}^{1 - \nu}$.
  Hence,
  \[
    \hat{C}_{k + 1}^{\nu / (1 - \nu)}
    [A_{k + 1}^{1 - \nu} - A_k^{1 - \nu}]
    \geq
    M^{2 \nu / (1 - \nu)} A_{k + 1}^{\nu}
    [A_{k + 1}^{1 - \nu} - A_k^{1 - \nu}].
  \]
  Thus, it suffices to show that
  \[
    A_{k + 1}^{\nu}
    [A_{k + 1}^{1 - \nu} - A_k^{1 - \nu}]
    \geq
    (1 - \nu) \alpha_{k + 1}.
  \]
  But this is indeed true, as for any $0 \leq t_1 \leq t_2$, by the concavity of
  $t \mapsto t^{1 - \nu}$, we have
  $t_2^{\nu} (t_2^{1 - \nu} - t_1^{1 - \nu}) \geq (1 - \nu) (t_2 - t_1)$,
  while $A_{k + 1} - A_k = \alpha_{k + 1}$ according to
  \cref{\LocalName{eq:UpperBoundOnRate}}.
\end{proof}

%% file: Appendix/AdditionalRelatedWork.tex
\section{Additional Related Work}
\label{sec:AdditionalRelatedWork}

Within the context of Problem~\eqref{eq:MainProblem},
we most commonly consider that $f$ and $\psi$ are both convex and $\psi$ is a
\emph{simple}, non-smooth function, such that we could solve
\cref{eq:MainProblem} efficiently by means of a \emph{proximity} function.

\paragraph{Classical methods:}
Focusing on the setting where $\psi$ is the indicator function of a compact set
$Q$, we can solve the problem at a rate of $O(1/\sqrt{k})$ when $f$ is
non-smooth while we can accelerate the convergence to $O(1/k^2)$ when $f$ has
Lipschitz continous gradients, i.e., $f$ is smooth.
These rates are shown to be tight when the gradient feedback is
noiseless~\parencite{nemirovskii1983problem}.
When the first-order oracle is stochastic with variance~$\sigma^2$, the lower
bounds imply a convergence rate
of~$O(\sigma / \sqrt{k})$~\parencite{nemirovski2009robust,lan2012optimal}.

The simple (sub-)gradient descent (GD)~\parencite{cauchy1847methode}, i.e.,
$x_{k+1} = x_k - \gamma_k g(x_k)$, with a sufficiently small $\gamma_0$ that
decays as $O(1/\sqrt{k})$ achieves $O(1 / \sqrt{k})$ rate for non-smooth
minimization.
Although this rate matches the information theoretic lower bounds, the same
method converges at an $O(1/k)$ rate under smoothness, which is suboptimal.
\textcite{nesterov1983method} introduced the idea of ``momentum'' and
proposed the first order-optimal algorithm, accelerated gradient descent (AGD),
which manages to decrease the error at a rate of $O(1/k^2)$.
Since then, various different interpretations of Nesterov's acceleration has
been proposed.
For a broad review of acceleration mechanisms, we refer the reader
to~\textcite{%
  nesterov2005smooth,%
  xiao2010dual,%
  tseng2008accelerated,%
  beck2009fast,%
  diakonikolas2018accelerated,%
  wang2018acceleration%
} and references therein.

An integral components of the classical methods, such as GD, AGD and its
variants, is the dependence on the knowledge of problem parameters, specifically
the Lipschitz constant of the problem.
When the step-size is not selected sufficiently small with respect to the
Lipschitz constant, then these methods are destined to diverge.
Similar arguments hold for stochastic methods such that the initial step-size
needs to be sufficiently small to guarantee convergence for smooth
problems~\parencite{nemirovskii1983problem}.
Additionally, the step-size must decay optimally at the rate
of~$O(1/ \sqrt{k})$, irrespective of the smoothness of the problem, to control
the effect of noise and ensure convergence to the set of
solutions~\parencite{robbins1951stochastic}.

\paragraph{Line-search methods:}
A fundamental technique to overcome the dependence on problem parameters is the
line-search machinery~\parencite{%
  armijo1966minimization,%
  wolfe1969convergence,%
  nocedal2006numerical%
},
which dynamically selects step-size every iteration by using local information.
There are several strategies such as exact line-search and backtracking
line-research, which could be implemented with appropriate ``sufficient
decrease'' and curvature conditions.
Essentially, line-search helps estimate a \emph{locally-valid} step-size,
enabling larger step-sizes than using the globally worst-case Lipschitz
constant.
When equipped with an appropriate line-search mechanism, GD and AGD could
achieve the same convergence rates without the need to know the Lipschitz
parameter.
However, this comes at the expense of an iterative search procedure which
demands function value evaluations per iteration of the line-search subroutine.
Similarly, stochastic variants of line-search are available, nonetheless, they
enforce extra assumptions on the objective and gradient
information~\parencite{paquette2020stochastic}.

%% file: Appendix/AdditionalExperiments/Main.tex
\section{Additional Experiments}

\begin{figure}[tb]
  \begin{center}
    \includegraphics[width=0.7\linewidth]{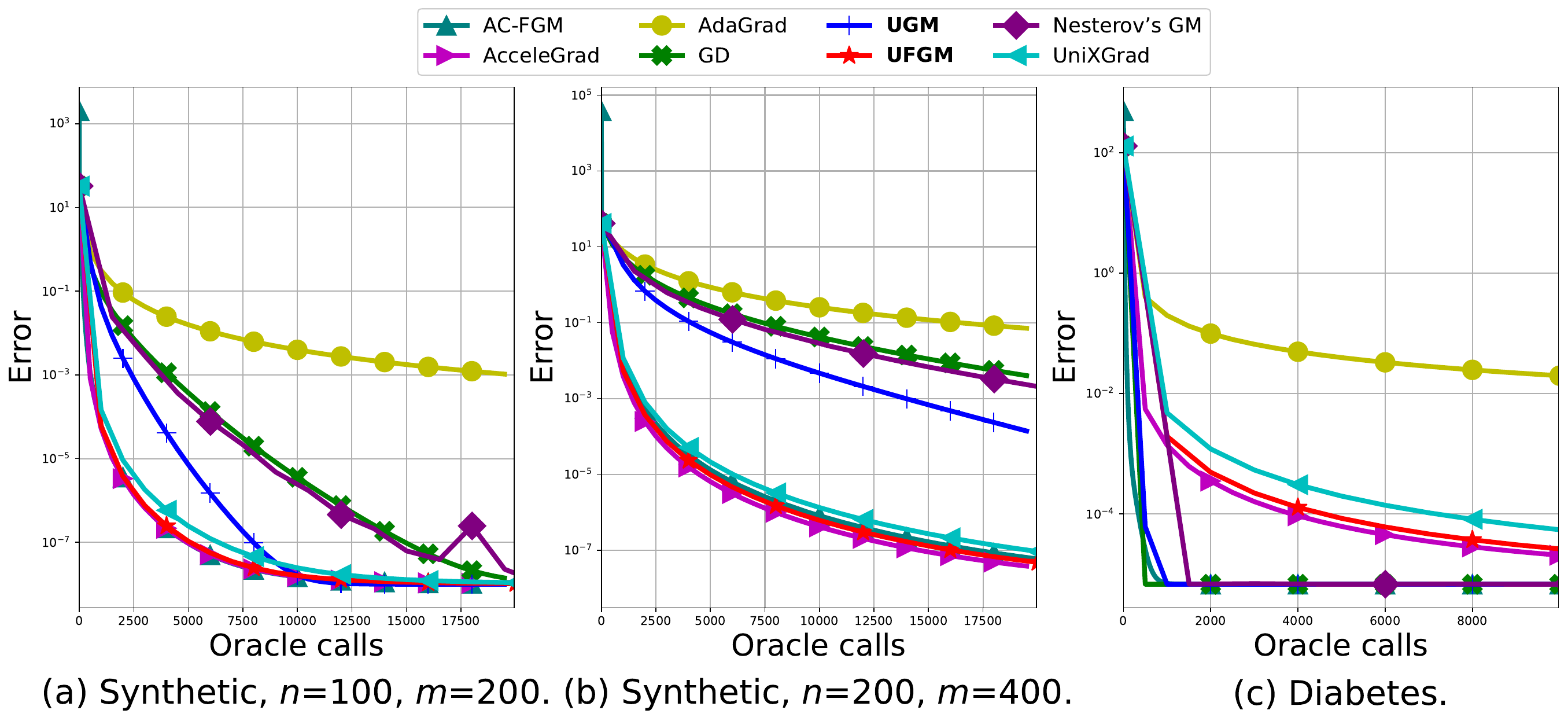}
  \end{center}
  \vspace{-6mm}
  \caption{%
    Comparison of different deterministic algorithms on convex optimization
    problems.
  }
  \label{fig:exp1determinstic}
\end{figure}
\begin{figure}[tb]
  \begin{center}
    \includegraphics[width=0.6\linewidth]{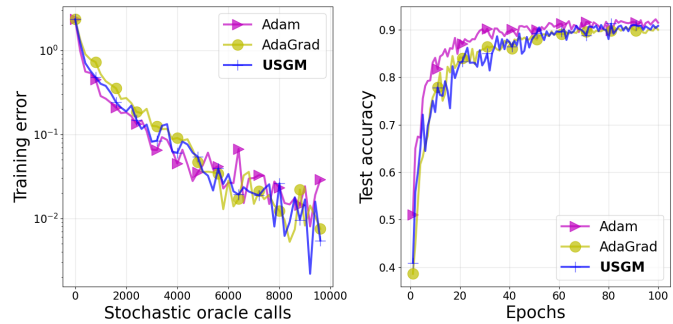}
  \end{center}
  \vspace{-7mm}
  \caption{%
    Comparison of different stochastic algorithms on non-convex optimization
    problems.
  }
  \label{fig:exp2resnet}
\end{figure}

In this section, we first elaborate on our experiments in the deterministic
setting.
We focus on the least-square problem in \cref{equ:ls}.
We first run the experiment on real-world diabetes dataset from LIBSVM.
Next, we consider a synthetic dataset, where we randomly generate an optimal
solution $x^*$ in the surface of a unit ball.
Next, we sample each element of $A$ from a uniform distribution
over~$\ClosedClosedInterval{0}{1}$ and set $b = A x^*$.
We test the proposed \cref{alg:GradientMethod}, denoted by UGM, and its
accelerated version (AUGM).
The baselines include GD, Nesterov's GM~\parencite{nesterov2015universal},
AdaGrad, UnixGrad~\parencite{kavis2019unixgrad},
AcceleGrad~\parencite{levy2018online} and AC-FGM~\parencite{li2023simple}.
In GD, we set the step size as $1/L$ while other methods are tuned via grid
search.
The result is presented in \cref{fig:exp1determinstic}, where we observe that
the proposed UGM shows better performance than UniXGrad and AcceleGrad.

Next, we include additional experiments on the stochastic setting.
To be specific, we train a ResNet18~\parencite{he2016deep} on
CIFAR-10~\parencite{krizhevsky2009learning}.
We select the mini-batch size as~$512$.
The step size of each method is tuned by a parameter sweep over
$\Set{10, 1, 0.1, 0.01, 0.001, 0.0001}$.
The diameter of the proposed method is tuned by sweeping over
$\Set{50, 35, 20, 10, 5}$.
We show the result in \cref{fig:exp2resnet}, where we can observe that the
proposed stochastic universal gradient method can be applied on non-convex
problems as well.